\documentclass[reqno,11pt]{amsart}
\usepackage{mathrsfs,graphicx,amsthm,amsfonts,amsxtra,xspace,array,cite,epic,mathdots,float}
\usepackage[left=2.5cm]{geometry}
\newcommand\version{public}
\newcommand\finalized{yes}
\usepackage{url}
\usepackage{enumerate}
\usepackage{asymptote}
\usepackage{amsmath,amssymb,amscd}

\usepackage{ifthen}     

\newcommand\pubpri[2]{%
\ifthenelse{\equal{\version}{public}}%
{{#1}}%
{\marginpar{\scshape\small Pubpri Alert}{#2}}}
\newcommand\pubprinoalert[2]{%
\ifthenelse{\equal{\version}{public}}%
{{#1}}%
{#2}}

\newcommand\ignore[1]{}

\usepackage{color}
\providecommand\wantcolor{yes}   %
\definecolor{backgroundyellow}{cmyk}{.2,.1,.8,.2}
\definecolor{backgroundblue}{rgb}{0,0,1}
\definecolor{backgroundred}{rgb}{1,0,0}
\definecolor{backgroundmagenta}{cmyk}{0,1,0,0}
\ifthenelse{\equal{\wantcolor}{no}}{\renewcommand\color[1]{}}{}

\newcommand\mysection{\section}

\newcommand\mysubsubsection[1]{%
		\subsubsection{\sffamily\upshape\mdseries #1}}

\newcommand\mysss{\mysubsubsection}

\newtheorem{annotation}{Annotation}
\newtheorem{theorem}[annotation]{
		Theorem}
\newtheorem{lemma}[annotation]{
		Lemma}
\newtheorem{definition}{
		Definition}
\newtheorem{corollary}[annotation]{
		Corollary}
\newtheorem{proposition}[annotation]{
		Proposition}
\newtheorem{example}{
		Example}
\newcommand\bexample{\begin{example}\begin{rm}}
\newcommand\eexample{\end{rm}\hfill$\Box$\end{example}}
\newtheorem{examplenobox}[annotation]{
		Example}
\newcommand\bexamplenobox{\begin{examplenobox}\begin{rm}}
\newcommand\eexamplenobox{\end{rm}\end{examplenobox}}
\newtheorem{exercise}[annotation]{
		Exercise}
\newcommand\bexercise{\begin{exercise}\begin{rm}}
\newcommand\eexercise{\end{rm}\end{exercise}}
\newtheorem{notation}[annotation]{
		Notation}
\newcommand\bnotation{\begin{notation}\begin{rm}}
\newcommand\enotation{\end{rm}\end{notation}}

\newtheorem{remark}[annotation]{
		Remark}
\newcommand\bremark{\begin{remark}
\begin{upshape}}
\newcommand\eremark{\end{upshape}
\end{remark}}
\newcommand\bdefn{\begin{definition}
\begin{upshape}}
\newcommand\edefn{\end{upshape}
\end{definition}}
\newtheorem{caveat}[annotation]{
		Caveat}
\newcommand\bcaveat{\begin{caveat}
\begin{upshape}}
\newcommand\ecaveat{\end{upshape}
\end{caveat}}
\newenvironment{caveatstar}{
\par\noindent{\scshape\bfseries
  Caveat: }\begin{rm}}{\end{rm}\newline} 
\newcommand\bcaveatstar{\begin{caveatstar}}
\newcommand\ecaveatstar{\end{caveatstar}}
\newenvironment{myproof}{%
\par\noindent{\scshape
  Proof: }\begin{rm}}{\hfill$\Box$\end{rm}\newline} 
\newcommand\bmyproof{\begin{myproof}}
\newcommand\emyproof{\end{myproof}}
\newenvironment{myproofnobox}{%
\par\noindent{\scshape Proof:}\begin{rm}}{\end{rm}\hfill\newline}
\newcommand\bmyproofnobox{\begin{myproofnobox}}
\newcommand\emyproofnobox{\end{myproofnobox}}
\newenvironment{solution}{%
\par\noindent{\scshape Solution: }\begin{rm}}{\hfill$\Box$\end{rm}\newline}
\newenvironment{solutionnobox}{%
\par\noindent{\scshape Solution: }\begin{rm}}{\end{rm}}
\newcommand\bsolution{\begin{solution}\begin{rm}}
\newcommand\esolution{\end{rm}\end{solution}}
\newcommand\bsolutionnobox{\begin{solutionnobox}\begin{rm}}
\newcommand\esolutionnobox{\end{rm}\end{solutionnobox}}
\newcommand\bthm{\begin{theorem}}
\newcommand\ethm{\end{theorem}}
\newcommand\bcor{\begin{corollary}}
\newcommand\ecor{\end{corollary}}
\newcommand\blemma{\begin{lemma}}
\newcommand\elemma{\end{lemma}}
\newcommand\bprop{\begin{proposition}}
\newcommand\eprop{\end{proposition}}
\newcommand\beqn{\begin{equation}}
\newcommand\eeqn{\end{equation}}
\newcommand\beqnstar{\begin{equation*}}
\newcommand\eeqnstar{\end{equation*}}

\newcommand\mtitle[1]%
{\marginpar{\sffamily\raggedright\upshape\small #1}}

\providecommand\finalized{no}
\ifthenelse{\equal{\finalized}{yes}}{}{\marginparwidth=2cm}
\ifthenelse{\equal{\finalized}{yes}}%
{\newcommand\mylabel[1]{\label{#1}}}%
{\newcommand\mylabel[1]{\label{#1}\marginpar{[{\ttfamily\upshape\tiny #1}]}}}
\ifthenelse{\equal{\finalized}{yes}}%
{\newcommand\checked[1]{}}%
{\newcommand\checked[1]{\marginpar{[{\ttfamily\upshape\tiny CHECKED: #1}]}}}
\ifthenelse{\equal{\finalized}{yes}}%
{\newcommand\spellchecked[1]{}}%
{\newcommand\spellchecked[1]{\marginpar{[{\ttfamily\upshape\tiny SPELLCHECKED: #1}]}}}

\providecommand\version{public}   
\ifthenelse{\equal{\version}{public}}%
{\newcommand\mcomment[1]{}}%
{\newcommand\mcomment[1]{\marginpar{{\raggedright\sffamily\upshape\small
\begin{spacing}{0.75} #1\end{spacing}}}}}
\ifthenelse{\equal{\version}{public}}%
{\newcommand\fcomment[1]{}}%
{\newcommand\fcomment[1]{\footnote{#1}}}
\ifthenelse{\equal{\version}{public}}%
{\newcommand\comment[1]{}}%
{\newcommand\comment[1]{{\small #1}}}

\newcommand\st{\,|\,}
\newcommand\fieldc{\mathbb{C}}
\newcommand\csastar{\csa^*}
\newcommand\csafstar{\csaf^*}
\newcommand\widehatp{\widehat{P}}
\newcommand\widehatpplus{\widehatp^+}
\newcommand\widehatw{\widehat{W}}
\newcommand\widehatwl{\widehatw_\Lambda}

\newcommand\germ{\mathfrak}

\newcommand{\be}{\begin{enumerate}}
\newcommand{\ee}{\end{enumerate}}
\newcommand{\beq}{\begin{equation}}
\newcommand{\eeq}{\end{equation}}
\newcommand{\complex}{\mathbb{C}}

\newcommand{\integers}{\mathbb{Z}}

\DeclareMathOperator{\wt}{wt}
\DeclareMathOperator{\GT}{GT}

\DeclareMathOperator{\tr}{trace}

\DeclareMathOperator{\supp}{supp} 
\DeclareMathOperator{\CL}{\mathfrak{c}} 
\DeclareMathOperator{\clb}{\mathcal{C}} 
\renewcommand\omega{\varpi}

\newcommand{\borel}{\widehat{\mathfrak{b}}}
\newcommand{\csa}{\mathfrak{h}}

\newcommand{\csaf}{\widehat{\mathfrak{h}}}
\newcommand\csahat{\widehat{\mathfrak{h}}}

\newcommand{\afw}{\widehat{W}} 
\newcommand{\afwex}{\widehat{W}_{\mathrm{ex}}} 
\newcommand{\trex}{T_{\mathrm{ex}}} 
\newcommand{\digautos}{\Sigma} 
\newcommand\sltwo{\mathfrak{sl}_2}
\newcommand{\currsl}[1][2]{\mathfrak{sl}_{#1}[t]}
\newcommand{\currsltwo}{\mathfrak{sl}_2[t]}
\newcommand{\afsl}{\widehat{\mathfrak{sl}_2}}
\newcommand{\sltwohat}{\widehat{\mathfrak{sl}_2}}
\newcommand\slm{\mathfrak{sl}_m}
\newcommand{\slmhat}{\widehat{\slm}}
\newcommand\lieghat{\widehat{\lieg}}

\newcommand{\scrc}{\mathfrak{P}}
\newcommand{\scrcstab}{\scrc^{\mathrm{stab}}}
\newcommand{\mcp}{\mathcal{Y}}
\newcommand{\stabasis}{\mathcal{B}}
\newcommand{\sign}{\epsilon}
\newcommand{\nbar}{{\overline{n}}}
\newcommand{\homh}{\mathfrak{t}} 
\newcommand{\xidag}{\xi^\dag}
\newcommand{\dprime}{{\prime\prime}}
\newcommand{\tsigma}{\tilde{\sigma}}
\newcommand{\tphi}{\tilde{\phi}}

\newcommand{\aform}[2]{\left(#1, #2\right)}
\newcommand{\bform}[2]{\left(#1 | #2\right)}
\newcommand{\cform}[2]{\langle #1, #2\rangle}
\newcommand\pogt{\scrc}
\newcommand\pogtl{\pogt(\lambda)}

\numberwithin{equation}{section} 
\begingroup\makeatletter\ifx\SetFigFont\undefined%
\gdef\SetFigFont#1#2#3#4#5{%
  \reset@font\fontsize{#1}{#2pt}%
  \fontfamily{#3}\fontseries{#4}\fontshape{#5}%
  \selectfont}%
\fi\endgroup%
\begin{document}
\title[Stability of CPL bases for $\currsltwo$]{Stability of the Chari-Pressley-Loktev bases\\ for local Weyl
  modules of $\currsltwo$}
\author{K.~N.~Raghavan, B.~Ravinder, Sankaran Viswanath}
\thanks{The second named author acknowledges support from CSIR under the SPM 
  Fellowship scheme.   The first and third named 
authors acknowledge support from DAE under a XII plan project}
\address{The Institute of Mathematical Sciences\\
CIT campus, Taramani\\
Chennai 600113, India}
\email{knr@imsc.res.in, bravinder@imsc.res.in, svis@imsc.res.in}

\begin{abstract}
We prove stability of the Chari-Pressley-Loktev bases for natural inclusions of local Weyl modules of the current algebra $\currsltwo$.       
These modules being known to be Demazure submodules in the level~$1$ representations of the affine Lie algebra $\sltwohat$,
we obtain, by passage to the direct limit,  bases for the level~$1$ representations themselves.
\end{abstract}

\maketitle
%
\section{Introduction}
\newcommand\lieg{\mathfrak{g}}
\noindent  
Let $\lieg$ be a finite-dimensional complex simple Lie algebra and $\lieg[t]$ the corresponding current algebra:   recall that the Lie bracket on $\lieg[t]$ is obtained from that on~$\lieg$ merely by the extension of scalars to the polynomial ring $\mathbb{C}[t]$ in one variable.
 Local Weyl modules, introduced by Chari and Pressley in \cite{CP}, are interesting finite-dimensional representations.   %
Corresponding to every dominant integral weight $\lambda$ of~$\lieg$,  there is one local Weyl module of~$\lieg[t]$ denoted by~$W(\lambda)$.  The $W(\lambda)$ is universal among finite-dimensional $\lieg[t]$-modules generated by a highest weight vector of weight~$\lambda$,  in the sense that any such module is uniquely  a quotient of~$W(\lambda)$~\cite{CP,CL}.

In~\cite{CP}, Chari and Pressley also produced nice monomial bases for local Weyl modules in the case $\lieg=\sltwo$.
Chari and Loktev in ~\cite{CL} clarified and extended the construction of these bases to the case $\lieg=\mathfrak{sl}_m$.  
Their work was motivated  by the following conjecture~\cite{CP,cplms} about the dimension of the local Weyl modules for $\lieg$ simply laced:  %
\begin{equation}\label{e:dimwl}
\textup{$\dim W(\lambda)=\prod_{i=1}^{\ell}\left(\dim W(\varpi_i)\right)^{a_i}$ \quad\quad\quad for $\lambda=a_1\varpi_1+\cdots+a_{\ell}\varpi_{\ell}$},
\end{equation} 
where $\ell$ is the rank and $\varpi_1$, \ldots, $\varpi_{\ell}$ the fundamental weights of~$\lieg$. %
Using their bases, Chari and Loktev were able to obtain, for $\lieg=\mathfrak{sl}_m$:
\begin{itemize}
\item the formula~(\ref{e:dimwl}) above,
\item an identification of~$W(\lambda)$ as a Demazure module of a level~$1$ representation of the untwisted affine Lie algebra $\lieghat$ corresponding to $\lieg$,
\item a fermionic formula for the graded character of $W(\lambda)$.
\end{itemize}

Later,   by very different means, %
Fourier and Littelmann in ~\cite{FL} obtained, for $\lieg$ simply laced, the identification of local Weyl modules as  level~$1$ Demazure modules, and then deduced formula \eqref{e:dimwl}.
Their methods do not however give explicit bases or fermionic character formulas. And for the base case 
$\lieg=\sltwo$ they have to refer to ~\cite{CP}.

In this paper, we investigate further the Chari-Pressley  bases for local Weyl modules of $\lieg=\sltwo$, taking into account the perspective gained from \cite{CL}.    The dominant integral weights for $\lieg=\sltwo$ being parametrized by the non-negative integers,  there is one local Weyl module~$W(n)$ for every integer $n\geq0$. 
Let us restrict ourselves in this introduction, for the sake of simplicity, to the case when $n$ is even.
The local Weyl modules then get identified with Demazure modules of the basic representation $L(\Lambda_0)$ of $\sltwohat$.    As such,  they are related by a chain of inclusions:  $W(0)\hookrightarrow W(2)\hookrightarrow W(4)\hookrightarrow\cdots$.  
It is natural to ask if the Chari-Pressley bases for the individual $W(n)$ respect these inclusions.    %
This question is the main focus of this paper.

As a first step,   we define---see Equation~\eqref{cldef}---a normalized version of the Chari-Pressley bases by replacing the powers in the monomials by divided powers and introducing a sign factor.   These normalized bases,  which we refer to throughout as the CPL (short for Chari-Pressley-Loktev) bases,   have better properties with respect to inclusions of local Weyl modules.   Indeed in our main result (Theorem~\ref{t:main}) we show that the CPL bases respect inclusions ``in the limit''.    Note that it is too much to expect any result of this nature without passage to the limit (Example~\ref{x:counter}).  

To state a little more precisely what we do,   let $\pogt(n)$ denote the parametrizing set of the CPL basis for $W(n)$: the elements of $\pogt(n)$ are pairs $(P,\pi)$ where $P$ is a Gelfand-Tsetlin pattern for $\sltwo$ with bounds $n,0$ and $\pi$ is a partition whose Young diagram fits into an $(n-p) \times p$ box, where $n-2p$ is the weight of the pattern $P$.      We first define a weight preserving embedding of $\pogt(n)$ into $\pogt(n+2)$ for each $n$, thereby obtaining a chain $\pogt(0)\hookrightarrow\pogt(2)\hookrightarrow\pogt(4)\hookrightarrow\cdots$.    Given an element $\xi$ of~$\pogt(n)$,   let $\xi_k$ be its image in~$\pogt(n+2k)$ (where $k$ is a non-negative integer),  and let $\CL(\xi_k)$ be the corresponding CPL basis element.   Consider the sequence $\CL(\xi_k)$, $k=0,1,2,\ldots$,  of elements in~$L(\Lambda_0)$. 
Our main result (Theorem~\ref{t:main}) implies that this sequence stabilizes for large~$k$.   In fact,  it says that $\CL(\xi_k)$ equals the stable value as soon as $k$ is such that the weight space of~$W(n+2k)$ corresponding to the weight  of $\xi$  equals that of~$L(\Lambda_0)$.
Passing to the direct limit, we obtain a basis for $L(\Lambda_0)$ consisting of the stable CPL basis elements (see~\S\ref{ss:lbasis}).    Moreover,   we obtain an explicit description of the stable CPL basis in terms of elements of the Fock space of the homogeneous Heisenberg subalgebra of $\sltwohat$ (Equations~\eqref{e:cplfock1},~\eqref{secondpart}).

As to the generalization of our results to the case $\lieg=\slm$,   we now briefly describe the issues that crop up.   %
The parametrizing set~$\pogtl$ of the Chari-Loktev basis for a local Weyl module~$W(\lambda)$ has a neat combinatorial description (see~\cite{rrv} for details):  namely,  it is the set  of {\em partition overlaid Gelfand-Tsetlin patterns} with boundary row~$\lambda$.    Further a natural normalization of the basis (analogue of Equation~\eqref{cldef}) suggests itself.
Denoting by $\theta$ the highest root of~$\lieg$,   the identification of local Weyl modules as Demazure modules (of some fundamental representation of~$\slmhat$) gives us a natural chain of inclusions  $W(\lambda)\hookrightarrow W(\lambda+\theta)\hookrightarrow W(\lambda+2\theta)\hookrightarrow\cdots$.    Mirroring this,  we have on the other hand a chain of weight preserving embeddings %
$\pogt(\lambda)\hookrightarrow\pogt(\lambda+\theta)\hookrightarrow\pogt(\lambda+2\theta)\hookrightarrow\cdots$.   
Given an element $\xi$ of~$\pogtl$,   let $\xi_k$ denote its image in~$\pogt(\lambda+k\theta)$,  and  $\CL(\xi_k)$ the corresponding normalized basis element.   %
It is tempting to conjecture, based on the evidence of the present paper,  that the sequence
$\CL(\xi_k)$ stabilizes,  and further that $\CL(\xi_k)$ equals the stable value once the weight space of~$W(\lambda+k\theta)$ corresponding to the weight  of $\xi$  stabilizes.  However, generalizing our methods to $\lieg=\slm$ presents formidable technical difficulties, and we hope to address these in future work.

This paper is organized as follows: in \S\ref{s:notation}, we set up the notation and recall some fundamental facts concerning local Weyl modules; in \S\ref{s:result}, we give the definition of the CPL bases, the statement of our main theorem (Theorem~\ref{mainthm}), as well as an application to constructing bases for level 1 representations of $\sltwohat$. Also described in  \S\ref{s:result} is  another variant of the Chari-Pressley basis, which is more natural when one thinks of Weyl modules as Demazure modules. The last section is devoted to the proof of Theorem~\ref{mainthm} in stages; the main case is when $n$ is even and the Gelfand-Tsetlin pattern is of weight zero, and the proof of this occurs in \S\ref{s:mainproof}. We show how the other cases can be reduced to this one using the Frenkel-Kac translation operators (\S\ref{fksec}) and automorphisms of $\sltwohat$ (\S\ref{s:odd}).

\medskip
\noindent
{\em Acknowledgments:}
The authors thank Vyjayanthi Chari for introducing them to the subject and for many helpful discussions.

\mysection{Notation and Preliminaries}\mylabel{s:notation}\label{notn}
\subsection{The affine Lie algebra $\sltwohat$}
Let $\mathfrak{sl}_2$ be the Lie algebra of 2$\times$2 trace zero
matrices over the field~$\fieldc$ of complex numbers with {\em standard
  basis\/}
\[ h=\begin{pmatrix}
1 & 0\\ 0 & -1
\end{pmatrix}, \quad\quad
x=\begin{pmatrix}
0 & 1\\ 0 & 0
\end{pmatrix},\quad\quad y=\begin{pmatrix}
0 & 0\\ 1 & 0
\end{pmatrix}.\]
Let $\csa=\complex h$ be the {\em standard\/} Cartan subalgebra
and $(A,B)\mapsto\tr (AB)$ the normalized invariant bilinear form on $\sltwo$.

\smallskip
Let $\complex [t,t^{-1}]$ be the ring of Laurent polynomials in an indeterminate $t$.
Let $\afsl$ be the affine Lie algebra defined by
\begin{center}
$\afsl=\mathfrak{sl_2}\otimes \complex [t,t^{-1}]\oplus \complex c\oplus \complex d,$ 
\end{center}
where $c$ is central and the other Lie brackets are given by 
\begin{align}
\left[A t^{m},B t^{n}\right] &= \left[A,B\right]  t^{m+n} + m\,\delta_{m,-n}\aform{A}{B}c, \label{lieb}\\
[d,A  t^{m}]&=m\left(A  t^{m}\right),\label{dlieb}
\end{align}
 for all $A,B \in \mathfrak{sl}_2$ and integers $m$, $n$:  here, as
 throughout the paper,  $At^s$ is shorthand for~$A \otimes t^s$.
We let $\csaf=\complex h \oplus\complex c\oplus\complex d$,
and regard $\csa^*$ as a subspace of $\csaf^*$ by setting $\cform{\lambda}{c}=\cform{\lambda}{d}=0$ for $\lambda \in \csa^*$.

Let $\alpha_0$, $\alpha_1$ denote the simple roots of $\afsl$ and let
$\alpha_0^{\vee} = c-h$, $\alpha_1^{\vee}=h$ be the corresponding
coroots. Let  $e_i$, $f_i$ ($i=0,1$) denote the Chevalley generators of
$\afsl$; these are given by 
\[ e_1 = x,\quad\quad f_1 = y,\quad\quad e_0 = yt,\quad\quad f_0 = xt^{-1} .\]
We have
\[\textup{$\cform{\alpha_1}{h}=2$,\quad\quad$\cform{\alpha_1}{c}=0$,\quad\quad
$\cform{\alpha_1}{d}=0$\quad\quad  and\quad\quad
$\cform{\alpha_0}{h}=-2$,\quad\quad$\cform{\alpha_0}{c}=0$,\quad\quad$\cform{\alpha_0}{d}=1$}.\]
Let $\delta = \alpha_0 + \alpha_1$ denote the null root, 
$\widehat{Q} = \integers\alpha_0 + \integers \alpha_1$ the root lattice, and 
$\widehat{Q}^+$ the non-negative integer span of $\alpha_0$, $\alpha_1$.
The weight lattice (resp.\ the set of dominant weights) is defined by 
\[ \widehat{P} \; (\text{resp.\ } \widehat{P}^+) = \{\lambda \in
\csaf^*: \cform{\lambda}{\alpha_i^\vee} \in \integers \; (\text{resp.\
}  \integers_{\geq 0}),\, i=0,1\}. \]
We define $\Lambda_0 \in \widehat{P}^+$  by $\cform{\Lambda_0}{h}=0, \cform{\Lambda_0}{c}=1, \cform{\Lambda_0}{d}=0$.

\smallskip
The Weyl group $\afw$ of $\afsl$ is the subgroup of $GL(\csaf^*)$ generated by the simple reflections $s_0,s_1$.
These are defined by $s_i(\lambda)=\lambda-\cform{\lambda}{\alpha_i^{\vee}} \,\alpha_i$  for $\lambda \in \csaf^*$, and $i=0,1$.  
\mcomment{$\langle\cdot,\alpha_0^\vee\rangle=(\cdot|\delta-\alpha_1)$ and $(\cdot
  |\delta)=\langle\cdot ,c\rangle$.} 
There is a non-degenerate, symmetric, bilinear $\afw$-invariant form
$\bform{\cdot}{\cdot}$ on $\csaf^*$, given by requiring that $\complex
\alpha_1$ be orthogonal to $\complex\delta + \complex \Lambda_0$,
together with the relations $\bform{\alpha_1}{\alpha_1} =2,
\bform{\delta}{\delta} = \bform{\Lambda_0}{\Lambda_0} = 0,
\bform{\delta}{\Lambda_0}=1$.

\smallskip
Given $\alpha \in \csa^*$, we define $t_\alpha \in GL(\csaf^*)$ by
\beq\label{tdef}
t_{\alpha}(\lambda)=\lambda+ \bform{\lambda}{\delta} \alpha -\bform{\lambda}{\alpha} \delta - \frac{1}{2} \bform{\lambda}{\delta} \bform{\alpha}{\alpha} \delta \;\;\text{ for } \lambda \in \csaf^*.
\eeq
\pubpri{}{
We have
\beqn\label{e:talpha}
t_\alpha(\alpha_1)=\alpha_1-(\alpha|\alpha_1)\delta\quad\quad
t_\alpha(\Lambda_0)=\Lambda_0+\alpha-\frac{1}{2}(\alpha|\alpha)\delta\quad\quad
t_\alpha(\delta)=\delta\quad\quad
\eeqn
As can be readily verified,  $t_\alpha$ is the
exponential of the nilpotent (of order $2$) linear operator $\tau_\alpha$ on $\csafstar$
given by \beqn\label{e:taualpha}
\tau_\alpha(\lambda)=(\lambda|\delta)\alpha-(\lambda|\delta)\alpha
\eeqn
The associations $\alpha\mapsto\tau_\alpha$ and
$\alpha\mapsto t_\alpha$ are
$W$-equivariant: $\tau_{w\alpha}=w\tau_\alpha w^{-1}$ for $t_{w\alpha}=wt_\alpha w^{-1}$ for
$\alpha\in\csastar$ and $w\in W$.     
}         %

The {\em translation subgroup} $T$ of $\afw$ is defined by $T = \{t_{j\alpha_1}: j \in \integers\}$.
We have $\widehat{W} = W\ltimes T$, where $W = \{1, s_1\}$ is the
underlying finite Weyl group. \mcomment{\raggedright
$s_0=s_1t_{-\alpha_1}$;   $wt_\alpha w^{-1}=t_{w\alpha}$ for
$\alpha\in\csastar$ and $w\in W$.
}
Now let $\omega_1 = \alpha_1/2$; then $Q= \integers \alpha_1$ and $P = \integers \omega_1$ are the root and weight lattices of the underlying $\mathfrak{sl}_2$. We also let $P^+ = \integers_{\geq 0} \, \omega_1$ be the set of dominant weights of the underlying finite-type diagram.

The extended affine Weyl group $\afwex$ is the semi-direct product 
$$\afwex =W\ltimes \trex,$$
where $\trex = \{t_{j\omega_1}: j \in \integers\}$.
Now consider the element $\sigma = s_1 t_{-\omega_1} \in \afwex$. 
This induces the diagram automorphism of the Dynkin diagram of $\afsl$; we have  $$\sigma\alpha_0 = \alpha_1,\, \sigma \alpha_1 = \alpha_0,\, \sigma\rho=\rho.$$
Here, $\rho\in \csaf^*$ is the Weyl vector, defined by
$\cform{\rho}{\alpha^\vee_i} = 1$ for $i=0,1$ and
$\cform{\rho}{d}=0$\pubpri{.}{: explicitly,
  $\rho=\omega_1+2\Lambda_0$.   Indeed, from~(\ref{e:talpha}), the images under
  $t_{-\omega_1}$ of $\alpha_1$, $\Lambda_0$, and $\delta$ are given by:
\beqn\label{e:forsigma} t_{-\omega_1}(\alpha_1)=\alpha_1+\delta\quad \quad\quad
t_{-\omega_1}(\delta)=\delta\quad\quad\quad t_{-\omega_1}(\Lambda_0)=\Lambda_0-\omega_1-\delta/4\quad\eeqn}
\mcomment{$\rho=\omega_1+2\Lambda_0$}
We also have
$\afwex = \afw \rtimes \digautos$, where $\digautos = \{1, \sigma\}$ is the subgroup generated by $\sigma$.

\subsection{The basic representation $L(\Lambda_0)$ of $\sltwohat$}
Given $\Lambda \in \widehat{P}^{+},$ let $L(\Lambda)$ be the
irreducible \mcomment{should one say ``integrable''?}
$\afsl$-module with highest weight $\Lambda$. It is the cyclic $\afsl$-module generated by $v_\Lambda$, with defining relations
\begin{align}
hv_\Lambda &= \cform{\Lambda}{h} v_\Lambda \;\;\;\forall h \in \csaf \label{defr1},\\
e_i v_\Lambda &= 0 \;\;\; (i=0,1)\label{defr2},\\
f^{\cform{\Lambda}{\alpha^{\vee}_i}+1}_i \,v_\Lambda &=0 \;\;\; (i=0,1).\label{defr3}
\end{align}\label{lq}

\noindent
It has weight space decomposition $L(\Lambda) = \oplus_{\mu \in
  \csaf^*} L(\Lambda)_{\mu}$\pubpri{. }{, where 
  \[L(\Lambda)_\mu:=\{v\in L(\Lambda)\st Hv=\mu(H)v \textup{ for all
    $H\in\csahat$}\}\]}%
The $\mu$ for which 
$L(\Lambda)_{\mu} \neq 0$ are the {\em weights of } $L(\Lambda)$.
The module $L(\Lambda_0)$ is particularly well-understood; 
the following well-known proposition describes the weight spaces
 of $L(\Lambda_0)$ \cite{K}.
\begin{proposition}\label{wtsbasicrep}
\be
\item The set of weights of $L(\Lambda_0)$ is 
$\{t_{j\alpha_1}(\Lambda_0)- d\delta\mid j\in\integers, 
d\in\integers_{\geq0}\}$.
\item $\dim\left(L(\Lambda_0)_{t_{j\alpha_1}(\Lambda_0)- d\delta}\right)=p(d)$, the number of partitions of $d$.
\ee
\end{proposition}

We let $\Lambda_1 = \sigma \Lambda_0$. Then, $\Lambda_0, \Lambda_1$ are (a choice of) fundamental weights corresponding to the coroots $\alpha^{\vee}_0, \alpha^{\vee}_1$, i.e., $\cform{\Lambda_i}{\alpha^\vee_j} = \delta_{ij}$ for $i,j \in \{0,1\}$. We let 
$v_{\Lambda_i}$ denote a highest weight vector of $L(\Lambda_i)$ for $i=0,1$.

\subsection{The current algebra and its Weyl modules}
Let $\complex [t]$ be the polynomial ring in an indeterminate $t$.
The current algebra
$\mathfrak{sl_2}[t]$=$\mathfrak{sl_2}\otimes\complex [t]$ is a Lie
algebra with Lie bracket is obtained from that of $\sltwo$ by
extension of scalars to $\complex[t]$:   
$[At^m,Bt^n]=[A,B]t^{m+n}$ for all  $A$, $B$ in $\sltwo$ and
non-negative integers $m$, $n$.    As such,  it is a subalgebra of $\sltwohat$.
\bdefn 
(see~\cite[\S1.2.1]{CL})
Given $n \in \mathbb{Z}_{\geq0}$, the {\em local Weyl module\/} $W(n)$ is the cyclic  $\mathfrak{sl_2}[t]$-module with generator $w_n$ and relations:
\beqn\label{e:wmrelations}(x t^{s})\,w_n=0, \quad\quad (h
t^{s+1})\,w_n=0, \quad\quad h\,w_n=nw_n, \quad\quad y^{n+1}\,w_n=0
\quad\quad\textup{for all $s \geq 0$}.\eeqn  
\edefn

\subsection{Weyl modules as Demazure modules}
We recall that the {\em standard\/} Borel subalgebra of $\afsl$ is
$$\borel= \mathfrak{sl}_2\otimes t \complex [t]\,\oplus\complex x \oplus \csaf.$$
Let $w$ be in $\afw$ and $\Lambda$ in $\widehat{P}^{+}$.     The 
weight space $L(\Lambda)_{w\Lambda}$ of $L(\Lambda)$
has dimension one (since two weights that are Weyl group conjugates
have the same multiplicities).

Define $V_w(\Lambda):={\mathfrak{U}}\borel\,\left(L(\Lambda)_{w\Lambda}\right)$ where $\mathfrak{U}\borel$ denotes the universal enveloping algebra of $\borel$.
Then, $V_w(\Lambda)$ is a $\mathfrak{U}\borel$-submodule of $L(\Lambda)$, called the 
{\em Demazure module} of $L(\Lambda)$ associated to $w$. More
generally,  given an element $w$ of the extended affine Weyl group
$\afwex$, we write $w = u \tau$ with $u \in \afw, \tau \in \digautos$
and define, following \cite{FL},
the associated Demazure module by $V_{w}(\Lambda):=V_u\left(\tau(\Lambda)\right)$.

We will consider the modules $V_{t_\lambda}(\Lambda_0)$ for $\lambda \in P$. 
It is convenient to use the notation of \cite{FL} and set 
\[ D(1, \lambda) := V_{t_{-\lambda}}(\Lambda_0).\]
\pubpri{ Since $\digautos=\{1, \sigma\}$,
the $D(1,\lambda)$ are Demazure modules for $L(\Lambda_0)$ (when
$\lambda \in Q$) or $L(\Lambda_1)$ (when $\lambda \not\in Q$).}{Put $\lambda=n\omega_1$.   Then $n$ is even iff $\lambda\in
  Q$ iff $t_{-\lambda}\in\afw$.   In this case,
  \beqn\label{e:demaeven}\textup{$t_{-n\omega_1}\Lambda_0=\Lambda_0-n\omega_1-\frac{n^2}{4}\delta$}\quad\quad\quad\textup{and}\quad\quad\quad
  \textup{
  $D(1,n\omega_1)={\mathfrak{U}}(\borel)L(\Lambda_0)_{(\Lambda_0-n\omega_1-\frac{n^2}{4}\delta)}$}\eeqn
 If $n$ is odd,  then
 $t_{-n\omega_1}=t_{-(n-1)\omega_1}t_{-\omega_1}=
(s_1t_{(n-1)\omega_1}s_1) t_{-\omega_1}=
(s_1t_{(n-1)\omega_1})(s_1 t_{-\omega_1})=
(s_1t_{(n-1)\omega_1})\sigma$,   and we have
\begin{gather}\label{e:demaodd1}\textup{$\sigma\Lambda_0=\Lambda_1=\Lambda_0+\omega_1-\frac{\delta}{4}$
  \quad \quad\quad
$s_1t_{(n-1)\omega_1}\Lambda_1=t_{-n\omega_1}\Lambda_0=
\Lambda_0-n\omega_1-\frac{n^2}{4}\delta$}\\
\label{e:demaodd2}\textup{$D(1,n\omega_1)=L(\Lambda_1)_{\Lambda_0-n\omega_1-\frac{n^2}{4}\delta}$}
\end{gather}
}
Further, $D(1,\lambda)$ is $\currsl$-stable (not just $\borel$-stable)
if, and only if, $\lambda \in P^+$; cf. \cite{FL}. 

The following theorem identifies the $\currsl$-stable Demazure modules
with the Weyl modules of the current algebra:  %
\begin{theorem} \textup{(Chari-Loktev)}
The local Weyl module $W(n)$ is isomorphic to the Demazure module $ D(1,n\omega_1)$, as modules of the current algebra $\currsl$.
\end{theorem}
This theorem was proved in \cite{CL} for $\currsl[n]$ and was
generalized in \cite{FL} to current algebras of types $A, D, E$.
The isomorphism of this theorem maps the generator $w_n$  of $W(n)$ to a vector  of $L(\Lambda_\nbar)$, which we will also denote~$w_n$. 
Here $\nbar$ is $0$ if $n$ is even and $1$ if $n$ is odd. 
By \cite[Corollary 1.5.1]{CL} (see also \cite[Corollary 1]{FL}), the
weight $\gamma$ of the vector $w_n \in L(\Lambda_\nbar)$ is a Weyl
conjugate of $\Lambda_\nbar$. Further, we must have $\cform{\gamma}{h}
= n$\pubpri{.}{, since weights are obviously preserved.}
\pubpri{It follows from \eqref{tdef} that  
$\gamma = t_{n\alpha_1/2}(\Lambda_0)$ (respectively
$t_{(n-1)\alpha_1/2}(\Lambda_1)$) if $n$ is even (respectively, if $n$
is odd). }{
A routine calculation shows that $\gamma$ is $t_{n\omega_1}\Lambda_0$
or $t_{(n-1)\omega_1}\Lambda_1$ accordingly as $n$ is even or
odd. Indeed,  we have $\afw=\{1,s_1\}\ltimes\{t_{j\alpha_1}\st
j\in\mathbb{Z}\}$, and from Equations~(\ref{e:demaeven},~\ref{e:demaodd1})
\beqn\label{e:gamma}
s_1 t_{j\omega_1}\Lambda_0=s_1 t_{(j-1)\omega_1}\Lambda_1=\Lambda_0-j\omega_1-\frac{j^2}{4}\delta 
\quad\quad
t_{j\omega_1}\Lambda_0=t_{(j-1)\omega_1}\Lambda_1=
\Lambda_0+j\omega_1-\frac{j^2}{4}\delta 
\eeqn
}

Since the $\gamma$-weight space of $L(\Lambda_{\nbar})$ is
one-dimensional, this isomorphism identifying the Weyl module as a
Demazure module is unique up to scaling. We will fix the following choice of $w_n$ for the rest of the paper:
\beq\label{wnchoice}
w_n :=\begin{cases} 
\left(xt^{-\frac{n}{2}}\right)^{\left(\frac{n}{2}\right)}\,v_{\Lambda_0}
& \text{  if $n$ is even},\\
\left(xt^{-\frac{n+1}{2}}\right)^{\left(\frac{n-1}{2}\right)}\,v_{\Lambda_1}
& \text{  if $n$ is odd}.
\end{cases}
\eeq
Here we have used the ``divided power notation'':
$X^{(p)}:=X^p/p!$.
It is clear that $w_n$ has weight $\gamma$; the fact that $w_n \neq 0$ will follow from Proposition~\ref{offdiag}(1) for $n$ even, and from the arguments of \S\ref{oddiso} for $n$ odd.
We will henceforth identify $W(n)$ with $D\left(1,n\omega_1\right)$ by the isomorphism defined by this choice of $w_n$, and think of $W(n)$ as a subspace of $L(\Lambda_\nbar)$.

\subsection{Inclusions of Weyl modules}
Let $\Lambda\in \widehatpplus$ and
$\widehatwl:=\{w\in\widehatw\st w\Lambda=\Lambda\}$.  For
elements $w_1\leq w_2$ of $\widehatw/\widehatwl$,  where $\leq$
denotes the Bruhat order on $\widehatw/\widehatwl$,  the Demazure
module $V_{w_1}(\Lambda)$ is included in $V_{w_2}(\Lambda)$ (as
submodules of $L(\Lambda)$).      Specializing to our case, we have, for $n$ even,
\beqn\label{e:weylinclude}
W(n)=V_{t_{-n\omega_1}}(\Lambda_0)\subseteq V_{t_{-(n+2)\omega_1}}(\Lambda_0)=W(n+2),
\eeqn
since $t_{-n\omega_1}\leq t_{-(n+2)\omega_1}=s_1s_0t_{-n\omega_1}$. For $n$ odd, we have $W(n) = V_{t_{-(n-1)\omega_1}s_1}(\Lambda_1)$, since
$t_{-\omega_1} = s_1 \sigma$. A similar argument to the above establishes $W(n) \subset W(n+2)$ in this case as well.
We thus have the following chains of embeddings:
\begin{align}
W(0)\hookrightarrow W(2)\hookrightarrow...\hookrightarrow &W(2n)\hookrightarrow W(2n+2)\hookrightarrow...\hookrightarrow L(\Lambda_0). \label{embeven}\\
W(1)\hookrightarrow W(3)\hookrightarrow...\hookrightarrow W&(2n+1)\hookrightarrow W(2n+3)\hookrightarrow...\hookrightarrow L(\Lambda_1). \label{embodd}
\end{align}

\mysection{The main results}\mylabel{s:result}
\subsection{Bases for Weyl modules}\mylabel{ss:basewn}
We first recall some results of \cite{CP} (see also \cite{CL}) which give a basis for the local Weyl module $W(n)$. 
We begin by introducing some notation. Let $\mcp$ denote the  set of
all integer partitions. Elements of $\mcp$ are infinite sequences $\lambda = (\lambda_1, \lambda_2, \lambda_3, \cdots)$ of non-negative integers such that (i) $\lambda_i \geq \lambda_{i+1}$ for all $i \geq 1$ and (ii) $\lambda_j =0$ for all sufficiently large~$j$. We let $|\lambda| = \sum_{i} \lambda_i$, and write $\lambda \vdash r$ to mean $\lambda \in \mcp$ with $|\lambda|=r$.  Let $\supp \lambda = \min \{j \geq 0: \lambda_{j+1} =0\}$. 
Given non-negative integers $a,b$, let 
\beqn\mcp(a,b) :=\{\lambda \in \mcp: \lambda_1 \leq b \text{ and }
\supp \lambda \leq a\}.\eeqn
We identify partitions with {\em Young diagrams\/} in the standard way:
the Young diagram corresponding to a partition~$\lambda$ is also
denoted $\lambda$ and consists of an arrangement of square
boxes, all of the same size (the sides are of unit length), 
numbering $|\lambda|$ in all,  arranged left-and top-justified,
$\lambda_1$ on the first row,  $\lambda_2$ on the second row (which is
below the first row), and so on:
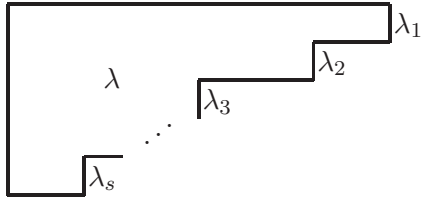
\begin{figure}[H]
\begin{center}
\setlength{\unitlength}{0.2in}
{\renewcommand{\dashlinestretch}{30}
\begin{picture}(6,5)(0,-5) 
\thicklines
\drawline(0,0)(10,0) 
\put(10.1,-.75){$\lambda_1$}
\drawline(10,0)(10,-1) 
\drawline(0,0)(0,-5) 
\drawline(0,-5)(2,-5) 
\put(2.1,-4.75){$\lambda_s$}
\drawline(2,-5)(2,-4) 
\drawline(2,-4)(3,-4) 
\drawline(8,-1)(10,-1) 
\drawline(8,-1)(8,-2) 
\put(8.1,-1.75){$\lambda_2$}
\drawline(5,-2)(8,-2) 
\drawline(5,-2)(5,-3) 
\put(5.1,-2.75){$\lambda_3$}
\put(3.5,-3.75){$\iddots$}
\put(2.5,-2.2){$\lambda$}
\end{picture}
}
\end{center}  
\end{figure}

\noindent
where $s=\supp \lambda$.
In this language, $\mcp(a,b)$ is the set of partitions whose Young diagrams fit
into a rectangular $a\times b$ box:
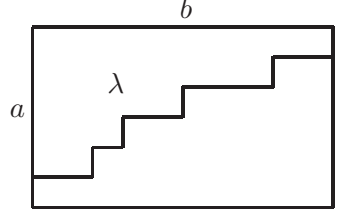
\begin{figure}[H]
\begin{center}
\setlength{\unitlength}{.4cm}
{\renewcommand{\dashlinestretch}{30}
\begin{picture}(6,5)(0,-5) 
\thicklines 
\drawline(0,0)(10,0) 
\put(4.9,0.3){$b$}
\drawline(10,0)(10,-6) 
\drawline(10,0)(10,-1) 
\drawline(0,0)(0,-6) 
\drawline(0,-5)(2,-5) 
\drawline(0,-6)(10,-6) 
\drawline(2,-5)(2,-4) 
\drawline(2,-4)(3,-4)
\drawline(3,-4)(3,-3)
\drawline(3,-3)(5,-3) 
\drawline(8,-1)(10,-1) 
\drawline(8,-1)(8,-2) 
\drawline(5,-2)(8,-2) 
\drawline(5,-2)(5,-3) 
\put(2.5,-2.2){$\lambda$}
\put(-0.75,-3){$a$}
\end{picture}
}
\end{center}
\end{figure}

Next, we define the set which will parametrize bases for local Weyl modules:
 \beq\label{scrcdef}
\scrc:=\{(m,k,\lambda): m, k \in \integers \text{ with } m \geq k \geq 0, \text{ and } \lambda \in \mcp(m-k,k)\}.
\eeq
In light of \cite{CL} %
a triple $(m,k,\lambda) \in \scrc$ should be thought of as the pair $(\GT_{m,k}, \,\lambda)$ where 
$$\GT_{m,k} = \begin{pmatrix}&k&\\m&&0 \end{pmatrix}$$ is  a Gelfand-Tsetlin pattern for $\mathfrak{sl}_2$.  
Associated to this pattern is a box of size $(m-k) \times (k-0)$, and the condition in \eqref{scrcdef} says that the Young diagram of $\lambda$ should fit into this box.

For each non-negative integer $n$, we also define $$\scrc(n) :=\{(m,k,\lambda) \in \scrc: m=n\}.$$
Given $\xi=(n,k,\lambda) \in \scrc(n)$ with $\lambda=(\lambda_1,\lambda_2, \lambda_3, \cdots)$, define the following element of $W(n)$:
\beq\label{bdef}
B(\xi) := \left(\prod_{i=1}^{n-k} yt^{\lambda_i}\right) w_n.
\eeq
We note that since $[yt^j,yt^k]=0$ for all $j,k \geq 0$, the order of terms in the product in equation \eqref{bdef} is irrelevant.
We now have the following important theorem due to Chari and Pressley
\cite{CP} (see also \cite{CL}): 
\begin{theorem}  \rm{(Chari-Pressley)}
Let $n\geq0$. Then $\{B(\xi)\mid \xi \in \scrc(n) \}$ is a basis for the local Weyl module $W(n)$.
\end{theorem}

\subsection{The CPL basis elements $\CL(\xi)$}
Our primary goal in this paper is to study the compatibility of the
Chari-Pressley bases for $W(n)$ with the chain of embeddings in equations \eqref{embeven} and \eqref{embodd}.
As a first step, we slightly modify the definition of these bases,
introducing normalization factors and parametrizing them by the
complements of partitions in $\mcp(n-k,k)$, rather than by the
partitions themselves. More precisely, given $\xi=(n,k,\lambda) \in
\scrc(n)$, define 
\beq \label{cldef}
\CL(\xi) := z(\xi) \left(\prod_{i=1}^{n-k} yt^{k-\lambda_i}\right) w_n.
\eeq
where $z(\xi)$ is a normalization factor. To specify $z(\xi)$, we first let $m_j :=\# \{i: \lambda_i =j\}$ denote the multiplicity of the part $j$ in $\lambda$ for each $j \geq 1$, and let $m_0 := n-k-\supp \lambda$. Then, we have $ \prod_{i=1}^{n-k} yt^{k-\lambda_i} = \prod_{j=0}^k \left( yt^{k-j}\right)^{m_j}$. 
The normalization factor is given by $$z(\xi) := \frac{(-1)^{[\frac{n}{4}]-[\frac{n-k}{2}]}}{\prod_{j=0}^k \; m_j !}.$$
Here, $[x]$ denotes the greatest integer less than or equal to $x$. We may also rewrite \eqref{cldef} in terms of divided powers; we have
$$\CL(\xi)=\sign(\xi) \,y^{(m_k)} \,(yt^1)^{(m_{k-1})} \cdots (yt^k)^{(m_0)} \,w_n,$$
where $\sign(\xi) = (-1)^{[\frac{n}{4}]-[\frac{n-k}{2}]}$.

Given $\xi = (n,k,\lambda) \in \scrc(n)$, with $s=\supp \lambda$, define $\lambda^{c} \in \mcp(n-k,k)$ by 
$$\lambda^{c} := (k,k,\cdots,k,k-\lambda_s, k-\lambda_{s-1}, \cdots,
k-\lambda_1,0,0,\cdots),$$ where the initial string of $k$'s is of
length $n-k-s$. The Young diagrams of $\lambda^{c}$ and 
$\lambda$, the latter rotated by $180^\circ$ and appropriately translated, are
complements of each other in the $(n-k) \times k$ box:
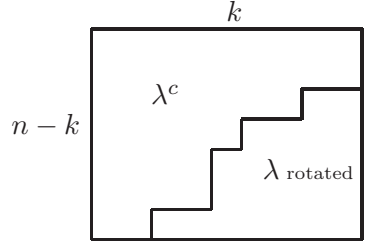
\begin{figure}[H]
\setlength{\unitlength}{0.4cm}
\begin{center}
{\renewcommand{\dashlinestretch}{30}
\begin{picture}(6,5)(0,-5.5)
\thicklines 
\drawline(0,0)(9,0)
\drawline(0,0)(0,-7)
\drawline(0,-7)(9,-7)
\drawline(9,0)(9,-7)
\put(4.5,.25){$k$}
\put(-2.65,-3.5){$n-k$}
\drawline(2,-7)(2,-6)(4,-6)(4,-4)(5,-4)(5,-3)(7,-3)(7,-2)(9,-2)
\put(2,-2.5){$\lambda^c$}
\put(5.7,-5){\textup{$\lambda${\tiny{ rotated}}}}
\end{picture}
}
\end{center}
\end{figure}

Letting $\xi^c = (n,k,\lambda^c)$, 
 it is clear that $\xi^c \in \scrc(n)$ and $\CL(\xi) = z(\xi) \, B(\xi^c)$. This of course implies that the set
$$\clb(n) :=\{\CL(\xi): \xi \in \scrc(n)\}$$
\smallskip
is also a basis for $W(n)$. We call this the {\em CPL basis} of $W(n)$.

We now view $W(n)$ as a subspace of $L(\Lambda_\nbar)$ as in equations \eqref{embeven} and \eqref{embodd}. The weight of $\CL(\xi)$ in $L(\Lambda_\nbar)$ is given by the following lemma.

\begin{lemma}\label{wtcl}
Let $\xi = (n,k,\lambda) \in \scrc$. Then
\be
\item
$ \text{Weight of } \CL(\xi) = t_{(k-n)\alpha_1}\left(\text{weight of } w_n\right) - |\lambda|\delta.$
\smallskip
\item 
If $n$ is even, the weight of $\CL(\xi)$  in  $L(\Lambda_0)$ is $ t_{(k-\frac{n}{2})\alpha_1} (\Lambda_0) - |\lambda|\delta $.
\smallskip
\item  If $n$ is odd, the weight of $\CL(\xi)$  in  $L(\Lambda_1)$ is $ t_{(k-\frac{n+1}{2})\alpha_1} (\Lambda_1) - |\lambda|\delta $.
\ee
\end{lemma}
\begin{proof}
From \eqref{cldef}, we have
\begin{align*}
\wt(\CL(\xi)) &= \wt(w_n) - (n-k)\alpha_1 + \delta \sum_{i=1}^{n-k}  (k-\lambda_i)\\
&= \wt(w_n) + (k-n)\alpha_1 + k(n-k)\delta - |\lambda|\delta.
\end{align*}
Let $\beta = \frac{n}{2}\alpha_1$ if $n$ is even, and
$\frac{n-1}{2}\alpha_1$ if $n$ is odd. Then $\wt(w_n) =
t_\beta(\Lambda_{\nbar})$. %
Since  $t_{(k-n)\alpha_1}$ and $t_\beta$ commute, the first part of the lemma is implied by the following identity, which can be verified directly using \eqref{tdef}:
$$ t_{(k-n)\alpha_1}(\Lambda_\nbar) = \Lambda_\nbar + (k-n) t_{-\beta}(\alpha_1) + k(n-k)\delta.$$
Assertions (2) and (3) are obvious from (1).
\end{proof}

\subsection{The main theorem: stability of the CPL bases}
We wish to study the compatibility of the bases $\clb(n)$  and $\clb(n+2)$ with respect to the embedding $W(n) \hookrightarrow W(n+2)$. As a first step, we define a weight preserving embedding at the level of the parametrizing sets of these bases. Define the map 
$\psi: \scrc \to \scrc$ by 
$$\psi(n,k,\lambda) = (n+2, k+1, \lambda).$$
This is well defined, since $\mcp(n-k,k)$ is a subset of $\mcp(n-k+1, k+1)$.
Further, $\psi$ is injective, and maps $\scrc(n)$ to $\scrc(n+2)$ for all $n$. Now, the following is immediate from Lemma~\ref{wtcl}. 
\begin{lemma}
Let $\xi \in \scrc(n)$. Then the basis vectors $\CL(\xi) \in W(n)$ and $\CL(\psi(\xi)) \in W(n+2)$ lie in the same weight space of $L(\Lambda_\nbar)$.
\end{lemma}

However, it is not true in general that $\CL(\xi)$ and $\CL(\psi(\xi))$ are equal as elements of  $L(\Lambda_\nbar)$, as the following example shows.
\bexample\mylabel{x:counter}
Let $\lambda$ be the partition $2+1$, i.e., $\lambda = (2,1,0,0,\cdots)$. Let $\xi = (4,2,\lambda)$. Then $\xi \in \scrc(4)$, and $\psi(\xi) = (6,3,\lambda)$. Using \eqref{cldef}, \eqref{wnchoice} and the commutation relations in $\afsl$, it is easy to compute:
\begin{align*}
\CL(\xi) &= \frac{1}{3} \left( ht^{-3} - (ht^{-1})^3\right) \, v_{\Lambda_0},\\
\CL(\psi(\xi)) &= \left(  ht^{-3} + ht^{-2}ht^{-1}\right) \, v_{\Lambda_0}.
\end{align*}
Both these vectors have weight $\Lambda_0 - 3\delta$. It is well known
that the vectors $ht^{-3}\,v_{\Lambda_0}$, 
\,$ht^{-2}ht^{-1}\,v_{\Lambda_0}$, $(ht^{-1})^3\,v_{\Lambda_0}$ form
a basis for the weight space $L(\Lambda_0)_{\Lambda_0 - 3 \delta}$. %
Thus, we conclude  $\CL(\xi)  \neq \CL(\psi(\xi))$.  %
\eexample

We will however see below that $\CL(\xi) = \CL(\psi(\xi))$ for all {\em stable} $\xi$. More precisely, let
\begin{equation}\label{fstabdef} %
\scrcstab(n) := \begin{cases}
\;\{(n, k, \lambda) \in \scrc(n): |\lambda| \leq \min(n-k, k)\} & \text{ if } n \text{ is even},
\vspace{2mm}\\
\;\{(n, k, \lambda) \in \scrc(n): |\lambda| \leq \min(n-k, k-1)\} & \text{ if } n \text{ is odd},
\end{cases}
\end{equation}
and $\scrcstab = \bigsqcup_{n \geq 0} \, \scrcstab(n)$.

\smallskip
We note that $\xi \in \scrcstab(n)$ implies $\psi(\xi) \in \scrcstab(n+2)$. The following is the main result of this paper.
\begin{theorem} \label{mainthm}\mylabel{t:main}
Let $n$ be a non-negative integer and $\xi =(n,k,\lambda)\in \scrcstab$. Then
$$\CL(\xi) = \CL(\psi(\xi)),$$
i.e., they are equal as elements of $L(\Lambda_\nbar)$.
\end{theorem}
This theorem is proved in \S\S\ref{s:mainproof}-\ref{s:odd}.

\subsection{Passage to the direct limit: a basis for~$L(\Lambda_0)$}\mylabel{ss:lbasis}
Theorem \ref{mainthm} allows us to construct a basis for $L(\Lambda_{p})$ ($p=0,1$) by taking the direct limit of the $\clb(n)$ (for $n \equiv p \pmod{2}$). We explain this below for $p=0$, the case $p=1$ being similar. Consider $L(\Lambda_0)$, and let $\mu = t_{j \alpha_1}(\Lambda_0) - d\delta \, (j \in \integers, d \in \integers_{\geq 0})$  be a weight of this module. Define
\beq\label{fmj}
\scrc_\mu := \{(n,k,\lambda) \in \scrc: k - \frac{n}{2} = j \text{ and } |\lambda| = d\}.
\eeq
We note that $\xi = (n,k,\lambda) \in \scrc_\mu$ forces $n$ to be even; further, 
it is clear from Lemma~\ref{wtcl} that $\CL(\xi)$ has weight $\mu$ iff $\xi \in \scrc_\mu$.

\smallskip
Now, let $\scrc_\mu(n) = \scrc_\mu \cap \scrc(n)$. This set parametrizes the basis elements of $W(n)$ of weight $\mu$. By \eqref{fmj}, the cardinality of $\scrc_\mu(n)$ is the number of partitions of $d$ which fit into a $(\frac{n}{2} -j) \times (\frac{n}{2} + j)$ box. Thus, for large enough $n$, $\scrc_\mu(n)$ contains exactly $p(d)$ (the number of partitions of $d$) elements; in particular this implies that $\psi$ induces a bijection of the sets $\scrc_\mu(n)$ and $\scrc_\mu(n+2)$. Further, it is also clear that for large $n$,  
every $\xi \in \scrc_\mu(n)$ is stable. More precisely, we have
\beq \label{fmust}
|\scrc_\mu(n)| = p(d) \text{ and } \scrc_\mu(n) \subset \scrcstab \text{ for all even } n \geq 2\left(d+ |j|\right).
\eeq
Choosing any such $n$, say $n = 2\left(d+|j|\right)$, we define the following (linearly independent) subset of $L(\Lambda_0)_\mu$:
$$\stabasis_\mu := \{ \CL(\xi): \xi \in \scrc_\mu(n) \}.$$
By Theorem~\ref{mainthm} and the remarks above, this is independent of
the choice of $n$. Since by Proposition~\ref{wtsbasicrep}, the
dimension of $L(\Lambda_0)_\mu$ is also $p(d)$, we conclude that
$\stabasis_\mu$ is a basis for the weight space
$L(\Lambda_0)_\mu$. Finally, to obtain a basis for $L(\Lambda_0)$, we take the disjoint union over the weights of  $L(\Lambda_0)$:
$$\stabasis := \bigsqcup_{\mu} \stabasis_\mu.$$
We may view $\stabasis$ as a direct limit of the CPL 
bases $\clb(n)$ ($n$ even) for the Demazure modules (=local Weyl modules) of $L(\Lambda_0)$. %
\subsection{A variation on the theme}
We note that the generator $w_n$ of $W(n) = D(1,n\omega_1)$ is not a lowest weight vector of the Demazure module $D(1,n\omega_1)$; while the lowest weight in 
$D(1,n\omega_1)$ is $t_{-n\omega_1}(\Lambda_0)$, the weight of $w_n$ is in fact $t_{n\omega_1}(\Lambda_0)$. From the basis $B(\xi)$ of equation \eqref{bdef}, it is easy to construct a basis consisting of monomials in the raising operators of the current algebra acting on a lowest weight vector $v_n$ of the Demazure module. Given $\xi=(n,k,\lambda) \in \scrc(n)$, with $\lambda=(\lambda_1,\lambda_2, \lambda_3, \cdots)$, define the following element of $W(n)$:
\beq\label{bbardef}
\overline{B}(\xi) := \left(\prod_{i=1}^{n-k} xt^{\lambda_i}\right) v_n.
\eeq
We now have:
\begin{proposition}\mylabel{p:bbar}
The set $\{\overline{B}(\xi)\mid \xi \in \scrc(n) \}$ is a basis for the local Weyl module $W(n)$.
\end{proposition}
\noindent
The proof appears in \S\ref{s:odd}. This basis also admits a
normalized version %
which exhibits similar stabilization behavior as the CPL basis.

\section{Proof of the main theorem}
\subsection{The key special case}\label{s:mainproof}
In this subsection we prove Theorem~\ref{mainthm} in the special case that $\xi = (n,k,\lambda) \in \scrcstab$ with $n$ even and $k = n/2$. In this case, the weight of $\CL(\xi)$ in $L(\Lambda_0)$ is $\Lambda_0 - |\lambda| \delta$.
From equations \eqref{cldef} and \eqref{wnchoice}, we have 
\beq \label{clxyform}
\CL(\xi) := z(\xi) \left(\prod_{i=1}^{k} yt^{k-\lambda_i}\right) \left(xt^{-k}\right)^{(k)}\,v_{\Lambda_0}.
\eeq

Now, let $\homh = \oplus_{n \in \integers} \complex  \, ht^n \oplus \complex  \,c$ denote the homogeneous Heisenberg subalgebra of $\afsl$. Recall that the subspace $\oplus_{p \geq 0} \,L(\Lambda_0)_{\Lambda_0 - p\delta}$ is invariant under $\homh$, and is isomorphic to the canonical commutation relations representation (Fock space) of $\homh$. Thus, each element of this subspace can be uniquely expressed as a polynomial in (the infinitely many variables) $ht^{-1}, ht^{-2}, \cdots$, acting on $v_{\Lambda_0}$ \cite{FK}. 
In particular, there is a unique polynomial $f_\xi(ht^{-1}, ht^{-2}, \cdots)$ such that
$$ \CL(\xi) = f_\xi(ht^{-1}, ht^{-2}, \cdots) \, v_{\Lambda_0}.$$

Our first goal is to determine $f_\xi$ explicitly by applying the {\em straightening rules} in $\mathfrak{U}\afsl$ to equation \eqref{clxyform}. We will then show that $f_\xi = f_{\psi(\xi)}$ for $\xi \in \scrcstab$, thereby establishing Theorem~\ref{mainthm} in this case.

\subsubsection{}
For $r \geq 1$, we let $[r]$:= $\{1,2,\cdots,r\}$.
Let $\pi \in \mcp$ be a partition such that $|\pi|=r$ and $\supp \pi =s$.
A {\em set partition} of $[r]$ of {\em type} $\pi$ is a collection $B = \{B_1, B_2, \cdots, B_s\}$ of pairwise 
disjoint subsets of $[r]$ such that $\cup_{i=1}^s B_i = [r]$ and $|B_i| = \pi_i$ for all $i \in [s]$.
We let $\mcp(\pi)$ denote the set of all set partitions of $[r]$ of type $\pi$.

Now, let $B=\{B_1, B_2, \cdots, B_s\} \in \mcp(\pi)$; given $\sigma \in S_r$ (the symmetric group on $r$ letters), $p=(p_1,p_2,\cdots,p_r) \in \mathbb{N}^r$ and $q=(q_1, q_2,\cdots,q_r) \in \mathbb{N}^r$,
define the following element of $\mathfrak{U} \homh$:
\beq \label{wdef}
W(B,\sigma; p, q):=\prod_{j=1}^{s}{ ht^{\sum_{i\in B_j}\left(p_i-q_{\sigma(i)}\right)}}.
\eeq

We also define 
\beq \label{hdef}
\mathcal{H}(\pi; p, q) := \frac{1}{\pi_1!...\pi_s!} \; \sum_{\substack{B\in \mcp(\pi) \\ \sigma\in S_r}} W(B,\sigma; p, q).
\eeq

With these notations we can state the following theorem.
\begin{theorem}\label{t1}
Let $r\geq 1$. For every triple $(p,q,v)$ with $p=(p_1,p_2,\cdots,p_r) \in \mathbb{N}^r$, $q=(q_1, q_2,\cdots,q_r) \in \mathbb{N}^r$ and $v \in L(\Lambda_0)$, satisfying
\be
\item $p_i< q_j$ for all $i,j \in [r]$,
\smallskip
\item $\sum_{i \in A} \, p_i \geq \sum_{j \in B} \, q_j$ for all subsets $A, B$ of $[r]$ such that $|A| = |B| + 1$,
\smallskip
\item $yt^{\left(\sum_{i \in A} \, p_i - \sum_{j \in B} \, q_j \right)} \,v=0$ for all subsets $A, B$ of $[r]$ such that $|A| = |B| + 1$,
\ee
we have
\beq\label{yteq}
\left(\prod_{i=1}^r yt^{p_i} \right)\left(\prod_{j=1}^r xt^{-q_j}\right) v =  
(-1)^{r}\sum_{\pi \vdash r}C(\pi)\,\mathcal{H}(\pi; p, q) \,v,
\eeq
where for $\pi= (\pi_1, \pi_2, \cdots)$, $C(\pi) = \prod_{i=1}^{\supp \pi}\, \pi_i!\,(\pi_i-1)!$.
\end{theorem}

\begin{proof}
We proceed by induction on $r$.
First, for $r=1$, consider $yt^{p_1}xt^{-q_1}\,v$. Since  $yt^{p_1} \, v=0$ and $p_1\neq q_1$, we have
$$ yt^{p_1}xt^{-q_1}\,v =[yt^{p_1}, xt^{-q_1}]\,v  =-ht^{{p_1}-{q_1}} \,v,$$
as required. Now let $r \geq 2$, and assume the result for $r-1$.
Consider $(\prod_{i=1}^r yt^{p_i})(\prod_{j=1}^r xt^{-q_j})\,v$. Since  $yt^{p_r} \, v=0$ and $p_r \neq q_j$ for all $j$, we may replace $yt^{p_r}\, (\prod_{j=1}^r xt^{-q_j})\,v$ by 
$$  [yt^{p_r}, \prod_{j=1}^r xt^{-q_j}]\,v =(-1)\sum_{l=1}^{r}\left(\prod_{j=l+1}^r xt^{-q_j}\right) ht^{p_r-q_l} \left(\prod_{j=1}^{l-1} xt^{-q_j}\right)\,v.$$
Next, using $[ht^{p_r-q_l}, xt^{-q_j}] = 2 xt^{-q_j-q_l+p_r}$, we can commute the $ht^{p_r-q_l}$ term past the 
$\left(\prod_{j=1}^{l-1} xt^{-q_j}\right)$. This yields
\beq\label{e2}
\begin{split}
(-1) \prod_{i=1}^r yt^{p_i}\,\prod_{j=1}^r xt^{-q_j}&\,v=
\sum_{l=1}^r \, \prod_{i=1}^{r-1} yt^{p_i} \prod_{\substack{j=1\\j \neq l}}^r xt^{-q_j} \left(ht^{p_r - q_l} \,v\right) \\
&+ 2 \sum_{\substack{l,m=1 \\ m< l}}^r \, \prod_{i=1}^{r-1} yt^{p_i} \prod_{\substack{j=1\\j \neq l,m}}^r xt^{-q_j}\,
(xt^{-q_m-q_l+p_r}) \, v.
\end{split}
\eeq

We now consider the first sum in equation \eqref{e2}. Fix $l \in [r]$ and let $p^\prime$ and $q^\prime$ denote the $r-1$ tuples obtained by deleting $p_r$ from $p$ and $q_l$ from $q$ respectively. We also let $v^\prime = ht^{p_r - q_l} \,v$. Then, we claim that the triple ($p^\prime, q^\prime, v^\prime)$ satisfies the hypotheses (1)-(3) of the theorem. The first two hypotheses are clear; now given $A \subset [r-1]$ and $B \subset [r]\backslash \{l\}$ with $|A| = |B| + 1$, we have
\beq
\begin{split}
yt^{\left(\sum_{i \in A} \, p_i - \sum_{j \in B} \, q_j \right)}v^\prime& = 
\left[yt^{\left(\sum_{i \in A} \, p_i - \sum_{j \in B} \, q_j \right)}, ht^{p_r - q_l}\right] \,v\\
 &=2yt^{\left(\sum_{i \in A\cup\{r\}} \, p_i - \sum_{j \in B\cup\{l\}} \, q_j \right)} \,v =0,
\end{split}
\eeq
thereby verifying hypothesis (3). By the induction hypothesis, we obtain 
\begin{equation}\label{e3}
\prod_{i=1}^{r-1} yt^{p_i} \prod_{\substack{j=1\\j \neq l}}^r xt^{-q_j} \left(ht^{p_r - q_l} \,v\right) = 
(-1)^{r-1}\sum_{\pi^\prime \vdash r-1}C(\pi^\prime)\,\mathcal{H}(\pi^\prime; p^\prime, q^\prime) \,ht^{p_r - q_l} \,v.
\end{equation}

The second sum in equation \eqref{e2} is treated analogously. Fix $l,m \in [r]$ with $m <l$ and let $q^\dprime$ denote the 
$r-1$ tuple obtained from $q$ by deleting $q_l, q_m$ and appending $q_l + q_m - p_r$. We also let $p^\dprime = (p_1, p_2, \cdots, p_{r-1})$ and $v^\dprime = v$. The triple $(p^\dprime, q^\dprime, v^\dprime)$ evidently satisfies the hypotheses of the theorem. Again, the induction hypothesis implies
\begin{equation}\label{e3b}
\prod_{i=1}^{r-1} yt^{p_i} \prod_{\substack{j=1\\j \neq l,m}}^r xt^{-q_j}\,
(xt^{-q_m-q_l+p_r}) \, v = (-1)^{r-1}\sum_{\pi^\dprime \vdash r-1}C(\pi^\dprime)\,\mathcal{H}(\pi^\dprime; p^\dprime, q^\dprime) \,v.
\eeq

Fix a partition $\pi \vdash r$, with $\pi = (\pi_1, \pi_2, \cdots)$ and $s = \supp \pi$. We can now find the coefficient 
$C(\pi)$ that occurs in equation \eqref{yteq}. 
Since the $yt^{p_i}$ commute pairwise and likewise the $xt^{-q_j}$, it is clear that the expression for
  $(\prod_{i=1}^r yt^{p_i})(\prod_{j=1}^r xt^{-q_j})\,v$ is invariant under the $S_r\times S_r$ action that permutes the 
$p_i$  and $-q_j$ among themselves. Thus, to find $C(\pi)$ it is enough to find the coefficient of the canonical word
\beq \label{canword}
ht^{\sum_{i=1}^{\pi_1}(p_i-q_i)}ht^{\sum_{i={\pi_1}+1}^{\pi_1+\pi_2}(p_i-q_i)} \cdots ht^{\sum_{i=\pi_1 + \cdots + \pi_{s-1}+1}^{r}
\;(p_i-q_i)}\,v
\eeq
in the RHS of \eqref{e2}.

We consider two cases (a) $\pi_s=1$, and (b) $\pi_s \geq 2$. In case (a), it is clear from equations 
\eqref{e2}, \eqref{e3} and \eqref{e3b} that the canonical word above  occurs only in 
$\prod_{i=1}^{r-1} yt^{p_i} \prod_{\substack{j=1}}^{r-1} xt^{-q_j} \left(ht^{p_r - q_r} \,v\right)$, and with 
coefficient $C(\pi^\prime)$ where $\pi^\prime=(\pi_1, \pi_2, \cdots, \pi_{s-1})\vdash {r-1}$.
Thus, 
\beq
C(\pi) = C(\pi^\prime) = \prod_{i=1}^{s-1}\pi_i!(\pi_i-1)!=\prod_{i=1}^{s}\pi_i!(\pi_i-1)!,
\eeq
since $\pi_s=1$.

In case (b), we have $\pi_s \geq 2$.
Again, examining equations \eqref{e2}, \eqref{e3} and \eqref{e3b}, it follows that the canonical word in this case occurs only in 
$$\prod_{i=1}^{r-1} yt^{p_i} \prod_{\substack{j=1\\j \neq l,m}}^r xt^{-q_j}\,
(xt^{-q_m-q_l+p_r}) \, v,$$
for all $l,m$ such that $$\pi_1 + \cdots + \pi_{s-1}+1 \leq m<l\leq r.$$ Each such pair $(l, m)$ contributes a 
 coefficient $C(\pi^\dprime)$ where $\pi^\dprime=(\pi_1, \pi_2, \cdots, \pi_{s-1}, \pi_s-1)\vdash {r-1}.$
Since $ r - \sum_{i=1}^{s-1} \pi_i = \pi_s$, we get
\begin{equation*}
C(\pi) = \binom{\pi_s}{2}2C(\pi^\dprime) = \pi_s(\pi_s-1)\left(\prod_{i=1}^{s-1}\pi_i!(\pi_i-1)!\right)(\pi_s-1)!(\pi_s-2)! 
= \prod_{i=1}^{s}\pi_i!(\pi_i-1)!,
\end{equation*}
as required. This proves Theorem \ref{t1}.

\end{proof}

\subsubsection{}
Let $\lambda=(\lambda_1, \lambda_2, \cdots)$ be a partition with $\supp \lambda = r \geq 1$.
Let $\pi \vdash r$ with $\supp \pi =s$, and let $B$=$\{B_1, B_2, \cdots, B_s\}$ be an element of $\mcp(\pi)$.
Define the following elements of $\mathfrak{U}(h\otimes t^{-1}\complex [t^{-1}])$:
\begin{align}
W(B,\lambda)&:=\prod_{p=1}^{s}{ ht^{-\sum_{j\in B_p}\lambda_j}}, \text{ and } \\
\mathcal{H}(\pi, \lambda)&:= \sum_{B\in{\mcp(\pi)}}W(B,\lambda). \label{hlamdef}
\end{align}

\bexample
$\mathcal{H}(\pi=(3), \lambda=(\lambda_1, \lambda_2, \lambda_3))= ht^{-(\lambda_1+\lambda_2+\lambda_3)}$.\\
$ \mathcal{H}(\pi=(2,1), \lambda=(\lambda_1, \lambda_2, \lambda_3))= ht^{-(\lambda_1+\lambda_2)}ht^{-\lambda_3}+
ht^{-(\lambda_1+\lambda_3)}ht^{-\lambda_2}+ht^{-(\lambda_2+\lambda_3)}ht^{-\lambda_1}$.\\
$\mathcal{H}(\pi=(1,1,1), \lambda=(\lambda_1, \lambda_2,
\lambda_3))= ht^{-\lambda_1}ht^{-\lambda_2}ht^{-\lambda_3}$. %
\eexample
\noindent

We now have the following important corollary to Theorem~\ref{t1}:
\begin{corollary}\label{c2}
Let $r\geq1$. Fix a partition $\lambda=(\lambda_1, \lambda_2, \cdots)$ with $\supp \lambda =r$. 
Then, for all $k\geq \mid\lambda\mid$, we have
\beq\label{c2eq}
\left(\prod_{i=1}^r yt^{k-\lambda_i}\right)\left(xt^{-k}\right)^{(r)}\,v_{\Lambda_0}= (-1)^{r}\,\sum_{\pi\vdash r}C^{\prime}(\pi)\,\mathcal{H}(\pi, \lambda)\,v_{\Lambda_0}.
\eeq
Here, for $\pi = (\pi_1, \pi_2, \cdots)$, $C^{\prime} (\pi)$ is given by $$C^{\prime} (\pi) = \prod_{i=1}^{\supp \pi}(\pi_i-1)!.$$

\end{corollary}
\begin{proof}
Consider $p = (p_1, p_2, \cdots, p_r)$ and $q = (q_1, q_2, \cdots, q_r)$ with $p_i=k-\lambda_i$ and $q_i = k$ for all $i\in [r]$. 
We claim that the triple $(p,q,v_{\Lambda_0})$ satisfies the hypotheses of Theorem \ref{t1}. To see this, observe first that $p_i < q_j$ for all $i,j \in [r]$. Further, 
if $A$ is a non-empty subset of $[r]$, we have
$$\sum_{i \in A} p_{i} = \sum_{i \in A} (k- \lambda_{i}) \geq |A|k -|\lambda| \geq (|A|-1)k.$$
Finally, the highest weight vector $v_{\Lambda_0} \in L(\Lambda_0)$ clearly satisfies 
$yt^{p}\,v_{\Lambda_0}=0$ $\forall$ $p\geq0$. Thus, by Theorem \ref{t1}, we obtain
\begin{equation}\label{e6}
\left(\prod_{i=1}^r yt^{k-\lambda_i}\right)\left(xt^{-k}\right)^{r}\,v_{\Lambda_0}= 
(-1)^{r}\sum_{\pi \vdash r}C(\pi)\mathcal{H}(\pi; p,q).
\end{equation}
with $C(\pi)= \prod_{i=1}^{\supp \pi}\pi_i!(\pi_i-1)!.$
Now since $q_j=k$ for all $j$, it is clear from equations \eqref{hdef} and \eqref{hlamdef} that 
\begin{equation}\label{e7}
\mathcal{H}(\pi; p,q) =\frac{r!}{{\pi_1}!{\pi_2}! \cdots {\pi_s}!}\mathcal{H}(\pi, \lambda).
\end{equation}
Equations \eqref{e6} and \eqref{e7} complete the proof.
\end{proof}
We observe that while the expression on the left hand side of equation \eqref{c2eq} depends on $k$, the one on the right hand side is independent of it. 
The fact that these two expressions are equal for $k \geq \mid \lambda \mid$ is precisely what leads to the stability properties of interest.

\subsubsection{}
We can now deduce the key special case of Theorem \ref{mainthm} that
we are after, namely  for $\xi$ of the form $(n,n/2,\lambda)$ with $n$ even.
Firstly, given a partition $\lambda \in \mcp$, let $r =\supp \lambda$ and  
$m_j(\lambda) =\# \{i: \lambda_i =j\}$ denote the multiplicity of the part $j$ in $\lambda$ for each $j \geq 1$. If $r \geq 1$, define
the following element of $\mathfrak{U}(h\otimes t^{-1}\complex [t^{-1}])$:
 \beq\label{fldef}
f_\lambda(ht^{-1}, ht^{-2}, \cdots):= \frac{(-1)^{r}}{\prod_{j \geq 1}\, m_j(\lambda)!} \, 
\sum_{\pi\vdash r} C^\prime(\pi)\mathcal{H}(\pi, \lambda),
\eeq
where $C^\prime(\pi) = \prod_{i=1}^{\supp \pi}(\pi_i-1)!$ as in Corollary \ref{c2}. If $r=0$, i.e., $\lambda$ is the empty partition, we let $f_\lambda:=1$.

Now, let $\xi = (n,\frac{n}{2},\lambda) \in \scrcstab$ with $n$ even.   As mentioned before, the weight of $\CL(\xi)$ in this case is $\Lambda_0 - |\lambda| \delta$. The expression of $\CL(\xi)$ as a polynomial in $ht^{-1}, ht^{-2}, \cdots$ acting on $v_{\Lambda_0}$
 is given by the following theorem. 
\begin{theorem}\label{t2}
Let $n$ be even and let $\xi = (n,n/2,\lambda) \in \scrcstab$.
Then
\begin{equation}\label{e:cplfock1}
\CL(\xi) = f_\lambda(ht^{-1}, ht^{-2}, \cdots) \, v_{\Lambda_0}.
\end{equation}
\end{theorem}
\begin{proof}
Let $r = \supp \lambda$ and $k=n/2$. If $r=0$, then $\CL(\xi)=(yt^{k})^{(n-k)}w_n=(yt^{k})^{(k)}(xt^{-k})^{(k)}\,v_{\Lambda_0}=v_{\Lambda_0}$, by Lemma~\ref{appendixlem}\,(1). Now, for $r\geq1$, 
$$(\prod_{i=1}^{r} yt^{k-\lambda_i})(yt^{k})^{(n-k-r)} \,w_n = (\prod_{i=1}^{r} yt^{k-\lambda_i})(yt^{k})^{(k-r)}(xt^{-k})^{(k)}\,v_{\Lambda_0}
= (\prod_{i=1}^{r} yt^{k-\lambda_i})(xt^{-k})^{(r)}\,v_{\Lambda_0},$$
again by Lemma~\ref{appendixlem}\,(1). The theorem now follows from this and equations \eqref{c2eq}, \eqref{cldef} and \eqref{fldef}.
\end{proof}
We now observe that $f_\lambda$ depends only on $\lambda$  and
not on $n$, thereby 
proving Theorem~\ref{mainthm} when $\xi$ is of the form $(n,n/2,\lambda)$:
\begin{corollary}\label{diagcase}
Let $n$ be even and let $\xi = (n,n/2,\lambda) \in \scrcstab$.
Then $ \CL(\xi) = \CL(\psi(\xi))$.
\end{corollary}

\subsection{The general case when $n$ is even}\label{fksec}
\noindent
We now turn to the remaining cases of Theorem~\ref{mainthm} for even $n$, i.e., $\xi = (n,k,\lambda) \in \scrcstab$ with $n$ even and $k \neq n/2$. We 
will now show how to reduce these to the case $k=n/2$ using the translation operators of Frenkel and Kac.
We recall the necessary facts from \cite{FK}, stated for our context.

 Let $\Delta:=\{\alpha_1,-\alpha_1\}$ be the set of all roots of $\mathfrak{sl_2}$, and set $E_{\alpha_1}:=x$ and $E_{-\alpha_1}:=y.$
 Let $(V,\pi)$ be an integrable representation of $\afsl$ with weight space decomposition $V = \oplus_{\mu \in \csaf^*} V_{\mu}$. For a real root $\alpha = \gamma + k \delta \;(\gamma \in \Delta, k \in \integers)$ of $\afsl$ we define
 \begin{equation}
 r_{\alpha}^{\pi}:=e^{-\pi(E_{\alpha})} e^{\pi(E_{-\alpha})} e^{-\pi(E_{\alpha})}.
 \end{equation}
where $E_{\alpha}:=E_{\gamma}t^{k}$. The operator  $r_{\alpha}^{\pi}$ is a linear automorphism of $V$ such that $r_{\alpha}^{\pi}(V_{\mu})=V_{s_{\alpha}(\mu)}$, where $s_\alpha \in 
\widehat{W}$ is the reflection defined by $\alpha$.
 
Next, we introduce the translation operators $T_\beta^\pi$ on $V$ for each $\beta \in Q = \integers \Delta$.
For $\gamma\in\Delta$, define
 \begin{equation}
 T_{\gamma}^{\pi}:=r_{\delta-\gamma}^{\pi}\,r_{\gamma}^{\pi}.
\end{equation}
and let $T_{p\gamma}^{\pi}:=(T_{\gamma}^{\pi})^{p}\, \forall \,p \in \mathbb{Z}_{\geq 0}$. 
These operators satisfy $T_\beta^\pi(V_\mu) = V_{t_\beta(\mu)}$ for all $\mu \in \csaf^*$, $\beta \in Q$.

We will only need these operators in two cases, namely when $(V, \pi)$ is either the adjoint representation 
or the basic representation of $\afsl$. We note that $T_\beta^{\mathrm{ad}}$ is in fact a 
Lie algebra automorphism of $\afsl$. For ease of notation, we will denote the translation operators corresponding to the basic representation simply by $T_\beta$, 
suppressing the $\pi$ in the superscript.

The key properties of the translation operators are given by
Propositions 1.2 and 2.3 of \cite{FK}. We
summarize them for our context below:

\begin{proposition} \label{fkprop}  {\rm (Frenkel-Kac)} 
\begin{enumerate}
\item $T_{p\alpha_1}^{\mathrm{ad}}(xt^{k})=xt^{k-2p}$ $\forall$ $p,k \in \mathbb{Z}$.\smallskip
\item $T_{p\alpha_1}^{\mathrm{ad}}(yt^{k})=yt^{k+2p}$ $\forall$ $p,k \in \mathbb{Z}$.\smallskip
\item $ T_{p\alpha_1} T_{q\alpha_1}=  T_{(p+q)\alpha_1}$ $\forall$ $p,q \in \mathbb{Z}$.\smallskip
\item $T_{p\alpha_1}AT_{-p\alpha_1}(v)=T_{p\alpha_1}^{\mathrm{ad}}(A)\,v$ $\forall$ $ A \in \afsl$, $v \in L(\Lambda_0)$, $p \in \mathbb{Z}$.\smallskip
\item $T_{p\alpha_1}(v_{\Lambda_0})=\prod_{i=1}^{p}xt^{-(2i-1)}\,v_{\Lambda_0}$ $\forall$ $p \geq 0$.\smallskip
\item $T_{p\alpha_1}(v_{\Lambda_0})=\prod_{i=1}^{-p}yt^{-(2i-1)}\,v_{\Lambda_0}$ $\forall$ $ p\leq 0.$
\end{enumerate}
\end{proposition}

The following is the key proposition that allows us to carry out a reduction to 
the case $k=n/2$.

\begin{proposition}\label{offdiag}
Let $n$ be even. Then, we have:
\be
\item $w_n=(-1)^{[\frac{n}{4}]} \,T_{n\alpha_1/2} (v_{\Lambda_0})$.
\item  Given $0 \leq k \leq n$, let $\gamma = (k -n/2)\alpha_1$. Then
\beq\label{firstpart}
(-1)^{[\frac{n}{4}]}\, w_n = (-1)^{[\frac{n-k}{2}]} \, T_\gamma(w_{2(n-k)}).
\eeq
\item Given $\xi = (n,k,\lambda) \in \scrcstab$, let $\xidag=(2(n-k), n-k, \lambda)$ and $\gamma(\xi) = (k-n/2)\alpha_1$. Then $\xidag \in \scrcstab$, and
\beq \label{secondpart}
\CL(\xi) = T_{\gamma(\xi)}\left(\CL(\xidag)\right).
\eeq
\ee
\end{proposition}

\begin{proof}
The proof of (1) will be given in the appendix (see Lemma~\ref{appendixlem}\eqref{lastpart}). 
Equation \eqref{firstpart} follows easily from (1) and Proposition~\ref{fkprop}\,(3).
To prove (3), we start with equation \eqref{cldef} and use Proposition~\ref{fkprop} again to obtain
\beq\label{finexp4}
T_{-\gamma(\xi)} \left(\CL(\xi)\right) = z(\xi) \left(\prod_{i=1}^{n-k} T_{-\gamma(\xi)}^{\mathrm{ad}}\left( yt^{k-\lambda_i}\right) \right) \left(T_{-\gamma(\xi)}w_n\right).
\eeq
Now, $T_{-\gamma(\xi)}^{\mathrm{ad}}\left( yt^{k-\lambda_i}\right) = yt^{n-k-\lambda_i}$. Further, it is clear from definition that $z(\xi) = (-1)^{[\frac{n}{4}]-[\frac{n-k}{2}]} z(\xidag)$. Plugging these and \eqref{firstpart} into \eqref{finexp4}, we obtain \eqref{secondpart}.
\end{proof}

We can now complete the proof of Theorem~\ref{mainthm} for $n$ even. Given $\xi = (n,k,\lambda) \in \scrcstab$, recall that $\psi(\xi) = (n+2, k+1, \lambda)$. It is now immediate from the definitions that
$$ \gamma(\xi) = \gamma(\psi(\xi)) \text{ and } \psi(\xidag) = \psi(\xi)^\dag.$$
Proposition \ref{offdiag} and Corollary~\ref{diagcase} now imply Theorem~\ref{t:main} for the case that $n$ is even.\qed

\subsection{The proof for odd $n$}\mylabel{s:odd}
In this section, we show how to reduce the case of $n$ odd to that of $n$ even,  using automorphisms of $\afsl$.
\subsubsection{}
Let $\tau$ be an automorphism of $\afsl$ such that $\tau \csaf = \csaf$. We have the induced action of $\tau$ on $\csaf^*$ by $\cform{\tau\lambda}{h} = \cform{\lambda}{\tau^{-1}h}$. 
Given an $\afsl$-module $V$, let $V^{\tau}$ denote the module with the twisted action $$g\circ v=\tau^{-1}(g)\,v \text{ for } g\in \afsl, v \in V.$$
Observe that for automorphisms $\tau_1, \tau_2$, we have $V^{\tau_1\tau_2} \simeq \left(V^{\tau_2}\right)^{\tau_1}$.

We now study the twisted actions on $L(\Lambda_0)$ by two specific automorphisms $\tsigma, \tphi$ of $\afsl$. First, recall from \S\ref{notn} that $\sigma=  s_1t_{-\omega_1} \in \widehat{W}_{ex}$ is an automorphism of the Dynkin diagram of $\afsl$; it swaps $\alpha_0, \alpha_1$ and fixes $\rho$. Consider the Lie algebra automorphism $\tsigma$ of $\afsl$ given by the relations
$$ \tsigma(e_i) = e_{1-i}, \; \tsigma(f_i) = f_{1-i}, \; \tsigma(\alpha^\vee_i) = \alpha^\vee_{1-i} \;(i=0,1) \text{ and } \tsigma(\rho^{\vee}) = \rho^\vee .$$
Here
 $\rho^\vee \in \csaf$ is the unique element for which  $\langle\alpha_0,{\rho}^{\vee}\rangle=1, \langle\alpha_1,{\rho}^{\vee}\rangle=1$ and $\langle\Lambda_0,\rho^{\vee}\rangle=0$. Clearly $\tsigma$ is an involution, and $$\tsigma(yt^m) = xt^{m-1},\, \tsigma(xt^m) = yt^{m+1}, \, \tsigma(ht^m) = -ht^m + \delta_{m,0} \,c \;\, \forall m \in \integers.$$
 Further, $\tsigma$ leaves $\csaf$ invariant, and its induced action on $\csaf^*$ coincides with $\sigma$.

To define the second automorphism $\tphi$, we employ the following simple lemma, which follows directly from  the Lie bracket relations \eqref{lieb}, \eqref{dlieb}.
\begin{lemma}\label{l}
Let $\phi$ be an automorphism of $\mathfrak{sl_2}$, which preserves the Killing form. Then $\phi$ can be extended to an automorphism $\tphi$ of 
$\afsl$ by defining $ \tphi(c)=c$, $\tphi(d)=d$ and $\tphi(At^{m})=\phi(A)\,t^{m}$ $\forall A \in\mathfrak{sl_2}, m\in \mathbb{Z}$. 
\end{lemma}

Now, consider the involution $\phi$ of $\mathfrak{sl}_2$ defined by 
\beq\label{phidef}
\phi(x)=y,\, \phi(y)=x, \,\phi(h)=-h.
\eeq
This preserves the Killing form, so by Lemma~\ref{l}, it extends to an automorphism (in fact, an involution) $\tphi$ of $\afsl$. 
It is again clear that (i) $\tphi$ preserves $\csaf$, and (ii) the induced action of $\tphi$ on $\csaf^*$ coincides with the simple reflection $s_1$.

\begin{proposition}\mylabel{l0l1}
With notation as above, we have (i) $L(\Lambda_0)^{\tsigma} \simeq L(\Lambda_1)$, and (ii) $L(\Lambda_0)^{\tphi} \simeq L(\Lambda_0)$.
\end{proposition}
\begin{proof}
To prove (i), consider the $\mathfrak{U}\afsl$-linear map $L(\Lambda_1) \to L(\Lambda_0)^{\tsigma}$ which sends $v_{\Lambda_1}$ to $v_{\Lambda_0}$. To show this is well defined, we only need to check that $v_{\Lambda_0} \in  L(\Lambda_0)^{\tsigma}$ satisfies the relations \eqref{defr1}-\eqref{defr3} for $\Lambda = \Lambda_1$. Since $\tsigma$ interchanges each pair $(e_0, e_1)$, $(f_0, f_1)$ and acts as $\sigma$ on $\csaf^*$, all three relations follow. Now, this map is a surjection, since $v_{\Lambda_0}$ generates $L(\Lambda_0)^{\tsigma}$. Since
$L(\Lambda_1)$ is irreducible, it must be an isomorphism.

A similar argument establishes (ii). We map $L(\Lambda_0) \to L(\Lambda_0)^{\tphi}$ by sending  $v_{\Lambda_0}$ to $v_{\Lambda_0}$. To show that this extends to a well-defined $\mathfrak{U}\afsl$-linear map on all of $L(\Lambda_0)$, we verify that $v_{\Lambda_0} \in L(\Lambda_0)^{\tphi}$ satisfies \eqref{defr1}-\eqref{defr3} for $\Lambda = \Lambda_0$. As above, \eqref{defr1} holds since the action of $\tphi$ on $\csaf^*$ coincides with $s_1$, and $s_1 \Lambda_0 = \Lambda_0$. Further, in $L(\Lambda_0)$, we have $\tphi^{-1}(e_0) v_{\Lambda_0} = xt \,v_{\Lambda_0} =0$ and $\tphi^{-1}(e_1) v_{\Lambda_0} = y v_{\Lambda_0} =0$. This establishes \eqref{defr2}. Finally, for \eqref{defr3}, we compute in $L(\Lambda_0)$:  $\tphi^{-1}(f_1) v_{\Lambda_0} = x v_{\Lambda_0} =0$, and $\tphi^{-1}(f_0)^2 v_{\Lambda_0} = \left(yt^{-1}\right)^2 v_{\Lambda_0}$. Since $yt^{-1}$ is in a real root space of $\afsl$, 
it is easy to see that this last term is also zero by a standard $\mathfrak{sl}_2$ argument (using the $\mathfrak{sl}_2$ spanned by $xt, yt^{-1}$ and $h + c$). The fact that it is an isomorphism follows as in (i).
\end{proof}

Let $\tau = \tsigma \tphi$. Then Proposition~\ref{l0l1} implies
$$ L(\Lambda_1) \simeq L(\Lambda_0)^{\tsigma} \simeq \left(L(\Lambda_0)^{\tphi}\right)^{\tsigma} \simeq  L(\Lambda_0)^{\tau}.$$ 
The isomorphism $F: L(\Lambda_1) \to L(\Lambda_0)^\tau$ maps $v_{\Lambda_1} \mapsto v_{\Lambda_0}$. It is then determined on all of $L(\Lambda_1)$ by $\afsl$-linearity, i.e., by the relation $$F(Xv) = \tau^{-1}(X)\,F(v) \quad \forall\ \  X \in \afsl,\ \  v \in L(\Lambda_1).$$

\subsubsection{} \label{oddiso}
We now prove Theorem \ref{mainthm} for $\xi=(n,k,\lambda)\in \scrcstab$ 
with $n$ odd. 
From \eqref{cldef} and \eqref{wnchoice}, we have
\begin{equation*}
\CL(\xi)= z(\xi)\left(\prod_{i=1}^{n-k} yt^{k-\lambda_i}\right) (xt^{-\frac{n+1}{2}})^{\left(\frac{n-1}{2}\right)} \,v_{\Lambda_1}.
\end{equation*} 
Applying the isomorphism $F$, we obtain
$$F\left(\CL(\xi)\right)= z(\xi)\left(\prod_{i=1}^{n-k}yt^{k-\lambda_i-1}\right)(xt^{-\frac{n-1}{2}})^{\left(\frac{n-1}{2}\right)} \,v_{\Lambda_0} = \CL(n-1,k-1,\lambda),$$
since $(-1)^{[\frac{n}{4}]}=(-1)^{[\frac{n-1}{4}]}$ for $n$ odd. Observe by \eqref{fstabdef} that $(n,k,\lambda) \in \scrcstab$ for $n$ odd, implies that 
$(n-1,k-1,\lambda)$ is also in $\scrcstab$.
Theorem \ref{mainthm} now follows for $\xi$ since we have already proved it for all even $n$. This completes the proof of that theorem in all cases. \qed
\subsubsection{}
Finally, we observe that the above ideas also give us a proof of Proposition~\ref{p:bbar}.   With notation as in that proposition, first let $n$ be even.   If  $G: L(\Lambda_0) \to L(\Lambda_0)^{\tphi}$ is the isomorphism constructed in the proof of Proposition~\ref{l0l1},  observe that $G(w_n)=(yt^{-\frac{n}{2}})^{\left(\frac{n}{2}\right)} \,v_{\Lambda_0}=v_n$, say, is a lowest weight vector of~$D(1,n\omega_1)$.  Further, for $\xi\in \scrc(n)$,   we have $G(B(\xi))=\overline{B}(\xi)$, thereby proving Proposition~\ref{p:bbar} in this case.  
The stable basis elements in this set up are simply the images of the $\CL(\xi)$, $\xi\in\scrcstab$, under the appropriate isomorphism $G$.
  The case of odd $n$ is analogous,  via the isomorphism $G':L(\Lambda_1)\to L(\Lambda_1)^{\tsigma\tphi\tsigma^{-1}}$.

\section*{Appendix}\noindent
The following lemma collects together the straightening rules in $L(\Lambda_0)$ that are used in the course of proving our main theorem. In principle, 
these can all be proved directly by working in the vertex operator realization of $L(\Lambda_0)$ \cite{FK}. The proofs below are simpler, and are included here for the sake of completeness.
\begin{lemma}\label{appendixlem}
Let $v_{\Lambda_0}$ denote a highest weight vector of $L(\Lambda_0)$. Then
\be
\item  $ (yt^{m})^{(l)} (xt^{-m})^{(m)} \, v_{\Lambda_0} =
  (xt^{-m})^{(m-l)}\,v_{\Lambda_0}$  $\forall \, 1\leqslant l \leqslant m$.
\smallskip
\item
$\prod_{i=1}^{r}xt^{2i-1}\,\prod_{i=1}^{r}yt^{-(2i-1)}\, v_{\Lambda_0}=v_{\Lambda_0} $ $\forall$ $r\in \mathbb{N}.$ 
\smallskip
\item $\prod_{i=1}^{r}yt^{2i-1}\,\prod_{i=1}^{r}xt^{-(2i-1)}\, v_{\Lambda_0}=v_{\Lambda_0} $ $\forall$ $r\in \mathbb{N}.$
\smallskip 
\item Let $p>q\geq0$ and let $v\in L(\Lambda_0)$ satisfy $yt^{p}\,v=ht^{p-q}\,v=0$. 
Then  $$yt^{p}\,(xt^{-q})^{(s)}\,v=-(xt^{-q})^{(s-2)}\,xt^{p-2q}\,v
\;\;\;\; \forall \, s\geq 2.$$
\item For $r\in 2 \mathbb{N}$ and $0\leq j \leq \frac{r}{2}$, we have
$$\left(\prod_{i=1}^{\frac{r}{2}+j}yt^{2i-1}\right)(xt^{-r})^{(2j)}\left(\prod_{i=1}^{\frac{r}{2}-j}xt^{-(2i-1)}\right)\,v_{\Lambda_0}=(-1)^{j}v_{\Lambda_0}.$$
\item For $r\in 2 \mathbb{N}-1$ and $0\leq j \leq \frac{r-1}{2}$, we have
$$\left(\prod_{i=1}^{\frac{r+1}{2}+j}yt^{2i-1}\right)(xt^{-r})^{(2j+1)}\left(\prod_{i=1}^{\frac{r-1}{2}-j}xt^{-(2i-1)}\right)\,v_{\Lambda_0}=(-1)^{j}v_{\Lambda_0}.$$ 
\item $\left(\prod_{i=1}^{r}yt^{2i-1}\right)(xt^{-r})^{(r)}\,v_{\Lambda_0}=(-1)^{[\frac{r}{2}]}\,v_{\Lambda_0}$ $\;\;\forall$ $r\in \mathbb{N}$.
\smallskip
\item\label{lastpart}
 $(xt^{-r})^{(r)}\,v_{\Lambda_0}=(-1)^{[\frac{r}{2}]} \, T_{r\alpha_1} (v_{\Lambda_0}) \neq 0$ $\;\;\;\forall$ $r\in \mathbb{N}$.
\ee
\end{lemma}
\smallskip
\begin{proof}
(1) Consider the Lie subalgebra of $\afsl$ spanned by $E := yt^m, F := xt^{-m}$ and $H:= -h + mc$. This is isomorphic to $\mathfrak{sl}_2$. Further, $E, F$ act locally nilpotently on $L(\Lambda_0)$, and we have $Hv_{\Lambda_0}  = m v_{\Lambda_0}$, $Ev_{\Lambda_0} =0$. The standard $\mathfrak{sl}_2$ calculation now shows $E^{(l)} F^{(m)} v_{\Lambda_0} = F^{(m-l)} v_{\Lambda_0}$.

\smallskip
\noindent
(2) Using Proposition \ref{fkprop},  it is easy to see that this is just a restatement of the identity $T_{-r\alpha_1} T_{r\alpha_1} v_{\Lambda_0} = v_{\Lambda_0}$.

\smallskip
\noindent
(3) As in (2), this is now the identity $T_{r\alpha_1} T_{-r\alpha_1} v_{\Lambda_0} = v_{\Lambda_0}$.

\smallskip
\noindent
(4) With the given hypotheses, we compute
\begin{equation*}
yt^{p}(xt^{-q})^{s}\,v = [yt^{p},(xt^{-q})^{s}]\,v 
= - \sum_{i=0}^{s-1}(xt^{-q})^{i}\,ht^{p-q}\,(xt^{-q})^{s-1-i}\,v.
\end{equation*}
We also have
$[ht^{p-q},(xt^{-q})^{u}] = 2u(xt^{-q})^{u-1}\,xt^{p-2q}$ for all $u
\geq 1$. Applying this to the above equation completes the proof.

\medskip
\noindent
(5) For $j=0$, this is just the statement of (3).
For $1 \leq j \leq \frac{r}{2}$, define  $v_j:=\prod_{i=1}^{\frac{r}{2}-j}xt^{-(2i-1)}\,v_{\Lambda_0}$. 
From weight considerations, it can be easily seen that $v_j$ satisfies $yt^{r+2j-1}\,v_j=0=ht^{2j-1}\,v_j$.
Thus, by (4), we obtain 
$$
yt^{r+2j-1}(xt^{-r})^{(2j)}\,v_j =-(xt^{-r})^{(2j-2)}xt^{-(r-2j+1)}\,v_j  =-(xt^{-r})^{(2j-2)} \, v_{j-1}.
$$
The result now follows by induction on $j$.

\smallskip
\noindent
(6) This is analogous to (5).

\smallskip
\noindent
(7) For $r$ even, put $j=\frac{r}{2}$ in (5) to obtain
\begin{equation*}
\prod_{i=1}^{r}yt^{2i-1}(xt^{-r})^{(r)}\,v_{\Lambda_0}=(-1)^{\frac{r}{2}}\,v_{\Lambda_0}.
\end{equation*}
Similarly, for $r$ odd, put $j=\frac{r-1}{2}$ in (6):
\begin{equation*}
\prod_{i=1}^{r}yt^{2i-1}(xt^{-r})^{(r)}\,v_{\Lambda_0}=(-1)^{\frac{r-1}{2}}\,v_{\Lambda_0}.
\end{equation*}
The equations above give us the desired result for all $r \in \mathbb{N}$.

\smallskip
\noindent
(8) Let $r\in \mathbb{N}$. Then 
$(xt^{-r})^{(r)}\,v_{\Lambda_0}$ and $ T_{r\alpha_1} (v_{\Lambda_0}) = \prod_{i=1}^{r} xt^{-(2i-1)}\,v_{\Lambda_0}$ belong to the 1-dimensional space  $L(\Lambda_0)_{\Lambda_0+r\alpha_1-r^{2}\delta}$, and so we must have
$$(xt^{-r})^{(r)}\,v_{\Lambda_0}=a\prod_{i=1}^{r} xt^{-(2i-1)}\,v_{\Lambda_0},$$ 
for some $a\in\complex$. But by (3) and (7), it follows that  $a=(-1)^{[\frac{r}{2}]}$ and that these vectors are non-zero.
\end{proof}


\bibliographystyle{bibsty-final-no-issn-isbn}
\addcontentsline{toc}{section}{References}
\ifthenelse{\equal{\finalized}{no}}{
\bibliography{abbrev,references}
}{
}

\end{document}